\def\VersionDateTime{29/January/2016, 17:02 GMT+9:00. Version $2.5$}
\documentclass[12pt]{amsart}
\usepackage{amssymb,amsthm,amsmath,amscd,rsfs}
\usepackage{latexsym}

\newcommand{\NN}{{\mathbb{N}}}
\newcommand{\ZZ}{{\mathbb{Z}}}
\newcommand{\QQ}{{\mathbb{Q}}}
\newcommand{\RR}{{\mathbb{R}}}

\newcommand{\OO}{{\mathcal{O}}}

\newcommand{\codim}{\operatorname{codim}}
\newcommand{\ch}{\operatorname{char}}

\newcommand{\Image}{\operatorname{Image}}

\newcommand{\Supp}{\operatorname{Supp}}

\newcommand{\Spec}{\operatorname{Spec}}

\newcommand{\cherncl}{{c}}

\newcommand{\Proof}{{\sl Proof.}\quad}
\newcommand{\QED}{{\unskip\nobreak\hfil\penalty50\quad\null\nobreak\hfil
{$\Box$}\parfillskip0pt\finalhyphendemerits0\par\medskip}}

\newcommand{\ord}{\operatorname{ord}}

\newcommand{\Tr}{\operatorname{Tr}}

\newcommand{\red}{\operatorname{red}}
\newcommand{\ndr}{\operatorname{nd-rk}}
\newcommand{\pr}{\operatorname{pr}}

\title
[Trace of abelian varieties 
and the
geometric Bogomolov conjecture]
{Trace of abelian varieties 
over function fields and the geometric Bogomolov conjecture
}
\author
{Kazuhiko Yamaki}
\date{\VersionDateTime}
\address
{Institute for Liberal Arts and Sciences,
Kyoto University, Kyoto, 606-8501, Japan}
\email{yamaki.kazuhiko.6r@kyoto-u.ac.jp}

\begin{document}

\theoremstyle{plain}
\newtheorem{Theorem}{Theorem}[section]
\newtheorem{Lemma}[Theorem]{Lemma}
\newtheorem{Proposition}[Theorem]{Proposition}
\newtheorem{Corollary}[Theorem]{Corollary}
\newtheorem{Main-Theorem}[Theorem]{Main Theorem}
\newtheorem{Theorem-Definition}[Theorem]{Theorem-Definition}
\theoremstyle{definition}
\newtheorem{Definition}[Theorem]{Definition}
\newtheorem{Remark}[Theorem]{Remark}
\newtheorem{Conjecture}[Theorem]{Conjecture}
\newtheorem{Claim}{Claim}
\newtheorem{Example}[Theorem]{Example}
\newtheorem{Key Fact}[Theorem]{Key Fact}
\newtheorem{ack}{Acknowledgments}       \renewcommand{\theack}{}
\newtheorem*{n-c}{Notation and convention}      
\newtheorem{citeTheorem}[Theorem]{Theorem}
\newtheorem{citeProposition}[Theorem]{Proposition}

\newtheorem{Step}{Step}

\renewcommand{\theTheorem}{\arabic{section}.\arabic{Theorem}}
\renewcommand{\theClaim}{\arabic{section}.\arabic{Theorem}.\arabic{Claim}}
\renewcommand{\theequation}{\arabic{section}.\arabic{Theorem}.\arabic{Claim}}

\def\Pf{\trivlist\item[\hskip\labelsep\textit{Proof.}]}
\def\endPf{\strut\hfill\framebox(6,6){}\endtrivlist}

\def\Pfo{\trivlist\item[\hskip\labelsep\textit{Proof of Proposition~\ref{ch-of-hyp}.}]}
\def\endPfo{\strut\hfill\framebox(6,6){}\endtrivlist}

\maketitle


\begin{abstract}
We prove
that the geometric Bogomolov conjecture
for any abelian varieties
is reduced to that for nowhere degenerate abelian varieties
with trivial trace.
In particular, 
the geometric Bogomolov conjecture
holds
for abelian varieties whose maximal nowhere degenerate abelian subvariety
is isogenous to a constant abelian variety.
To prove the results, we 
investigate closed subvarieties of 
abelian schemes over constant varieties,
where constant varieties are varieties over
a function field which can be defined over the constant
field of the function field.
\end{abstract}

\section{Introduction}

In this paper,
we contribute to the
geometric Bogomolov conjecture for abelian varieties
by
investigating closed subvarieties of 
 abelian schemes over constant varieties,
where constant varieties are varieties over
a function field which can be defined over the constant
field of the function field.

\subsection{Background}

The geometric Bogomolov conjecture 
for abelian varieties
is an analogue of the theorem established
by Ullmo \cite{ullmo} and Zhang \cite{zhang2},
called the (arithmetic) Bogomolov conjecture.
To describe this theorem,
let $K$ be a number field.
Let $A$ be an abelian variety over $\overline{K}$.
A line bundle on $A$ is said to be
\emph{even}
if
it is preserved under
the pull-back by the inverse $a \mapsto -a$ of $A$.
To an even ample line bundle on $A$,
one associates
the canonical height function,
also called the N\'eron--Tate height,
which is known 
to be a semi-positive definite quadratic form
on the additive group $A \left( \overline{K} \right)$.
For an $x \in A \left( \overline{K} \right)$,
the value at $x$ of the canonical height function is called the (canonical)
height of $x$.
We say that a closed subvariety $X$ of $A$
\emph{has dense small points} if,
for any $\epsilon > 0$,
the set of points of $X \left( \overline{K} \right)$
whose canonical 
heights are not greater than $\epsilon$
is dense in $X$.
It is known that the density of small points does not depend on the choice
of even ample line bundles on $A$.
Since the canonical height function
is a quadratic form,
it follows
that a torsion subvariety has dense small points,
where a torsion subvariety is the translate of an abelian subvariety
by a torsion point.
The (arithmetic) Bogomolov conjecture for $A$,
which is the theorem of Ullmo and Zhang,
claims that the converse also holds, that is,
any irreducible closed subvariety $X$ of $A$
with dense small points is a torsion subvariety.

The theorem of Ullmo and Zhang is 
generalized by Moriwaki in \cite{moriwaki5}
to the case where $K$ is a field finitely generated over $\QQ$.
He considered
the arithmetic heights
associated to ``big'' polarizations
(when the transcendence degree of $K$ over $\QQ$ is positive)
and established the Bogomolov conjecture with respect to the canonical
height arising from a big polarization.
We remark that this kind of heights
are different from the classical ``geometric'' heights over function fields;
they are
more
arithmetic.

It is then quite natural to ask
whether there is an analogue 
of the theorem of Ullmo and Zhang
over a function field with respect to the classical geometric height.
Now, 
let $K$ be the function field of a
(normal projective) variety over an algebraically closed field $k$.
In \cite{gubler2},
Gubler has proved that the same statement as the theorem
holds for abelian varieties over $\overline{K}$
which are totally degenerate at some place.
However,
this statement 
does not hold for \emph{any} abelian variety in fact.
For example, one sees that, if $A$ is a constant abelian variety, that is,
$A = \widetilde{A} \otimes_k \overline{K}$ for some abelian variety 
$\widetilde{A}$ over $k$,
then,
for any closed subvariety $\widetilde{X}$ of $\widetilde{A}$,
the subvariety $X := \widetilde{X} \otimes_k \overline{K}$
of $A$ has dense small points.
More generally,
for any abelian variety $A$ over $\overline{K}$,
let $\left( \widetilde{A}^{\overline{K}/k} , \Tr_{A} \right)$
be the $\overline{K}/k$-trace of $A$
(cf. \S~\ref{sect:NC}).
Then, one sees that, for any closed subvariety $\widetilde{Y}$ of 
$\widetilde{A}^{\overline{K}/k}$, the closed subvariety
$\Tr_A \left( \widetilde{Y} \otimes_k \overline{K} \right)$
of $A$ has dense small points.
This means that, over function fields,
an abelian variety may in general have 
an irreducible closed subvariety
that is not a torsion subvariety but has dense small points.
It should be remarked that,
in the case of Gubler, 
it follows from the assumption of total degeneracy
that
$A$ has trivial $\overline{K}/k$-trace.
Therefore, there are no
such 
subvarieties,
and actually, he 
succeeded
in
establishing
the same statement 
for such abelian varieties.

Inspired by the work
due to Gubler,
we formulate in \cite{yamaki5} a Bogomolov conjecture
for any abelian variety over a function field
generalizing Gubler's theorem,
which we
call the \emph{geometric Bogomolov conjecture
(for abelian varieties)}
(cf. Conjecture~\ref{GBCforAV}).
This conjecture
claims that
any
irreducible closed subvariety
with
dense small points
is
a
\emph{special subvariety},
which has been defined as the
sum of a torsion subvariety
and a closed subvariety coming from the constant subvarieties of the 
$\overline{K}/k$-trace.
In \cite{yamaki5}, we established
a
partial answer to the conjecture,
and this result
generalizes
the theorem of the totally degenerate case by Gubler
in fact.
Further development can be found in \cite{yamaki6},
and more details of it are resumed in 
\S~\ref{sect:GBC}.
We remark that the geometric Bogomolov conjecture is still an
open problem.

We put a few comments on another
version of the conjecture, the geometric
Bogomolov conjecture for curves.
It
is a restricted version of 
the 
conjecture for abelian varieties,
claiming that any projective curve of genus more than $1$
embedded in its jacobian variety does not have dense small points.
When $\ch (k) = 0$ 
and $K$ 
has
transcendence degree $1$ over $k$,
this is proved
by Cinkir in \cite{cinkir}\footnote{
In fact, its effective version is proved here.}.
In positive characteristics, 
although
one can find some partial answers in
\cite{zhang1, moriwaki0, moriwaki3, yamaki1, yamaki4},
it
has not yet been solved.

\subsection{Notation and convention} \label{sect:NC}

Let $k$ be an algebraically closed field,
let
$\mathfrak{B}$ be an irreducible normal projective variety 
of dimension $b \geq 1$ over $k$,
and let $\mathcal{H}$ be an ample line bundle on $\mathfrak{B}$.
Let
$K$ be the function field 
of 
$\mathfrak{B}$
and let $\overline{K}$ be an algebraic closure of $K$.
All of them are fixed throughout this article.
Any finite extension of $K$ will be
taken in $\overline{K}$.

Let $M_{K}$ be the set of points of $\mathfrak{B}$
of codimension $1$.
Each element of $M_K$ is called a \emph{place} of $K$.
For any $v \in M_{K}$, the local ring $\OO_{\mathfrak{B}, v}$ 
is a discrete valuation
ring with fractional field $K$,
and we let
$\ord_{v} : K^{\times} \to \ZZ$ denote the order 
function.
This
gives rise to
a
non-archimedean 
value $| \cdot |_{v , \mathcal{H}}$ on $K$
normalized
in such a way that
\addtocounter{Claim}{1}
\begin{align} \label{eq:norm}
| x |_{v , \mathcal{H}} := 
e^{- \deg_{\mathcal{H}} ( \overline{v} )
 \ord_{v} (x)}
\end{align}
for any $x \in K^{\times}$,
where 
$\deg_{\mathcal{H}} ( \overline{v} )$
denotes the degree with respect to $\mathcal{H}$
of the closure $\overline{v}$ of $v$
in $\mathfrak{B}$.
It is well known that the set 
$
\{ | \cdot |_{v , \mathcal{H}} \}_{v \in M_{K}}$ 
of values satisfies the product formula,
and hence 
the notion of (absolute logarithmic) heights
with respect to this set of values
is defined (cf. \cite[Chapter~3 \S~3]{lang2}).

For a finite extension $K'$ of $K$ in $\overline{K}$,
let $\mathfrak{B}_{K'}$ be the normalization of $\mathfrak{B}$ in $K'$
and
let $M_{K'}$ be the set of points of $\mathfrak{B}_{K'}$ of codimension $1$.
For a finite extension $K'' / K'$, we have a natural surjective map $M_{K''} \to M_{K'}$,
and thus we obtain a inverse system $(M_{K'})_{K'}$,
where $K'$ runs through the finite extensions of $K$ in $\overline{K}$.
Set $M_{\overline{K}} := \varprojlim_{K'} M_{K'}$.
We call an element of $M_{\overline{K}}$ a \emph{place of $\overline{K}$}
(cf. \cite[\S~6.1]{yamaki6}).
Each $v \in M_{\overline{K}}$ gives a unique
absolute value on
$\overline{K}$ which extends $| \cdot |_{v_{K} , \mathcal{H}}$,
where $v_K$ is the image of $v$ by the natural map $M_{\overline{K}}
\to M_{K}$.
We denote by $\overline{K}_{v}$ the completion 
of $\overline{K}$ with respect to that absolute value.

Let $F / k$ be a field extension
and let
$X$ be a scheme over $F$.
For a field extension $\mathfrak{F}/F$,
we write 
${X} \otimes_F \mathfrak{F} := 
{X} \times_{\Spec (F)} \Spec (\mathfrak{F})$.
For a morphism $\phi : {X} \to {Y}$ of schemes over $F$,
we write $\phi \otimes_F \mathfrak{F} : 
{X} \otimes_F \mathfrak{F} \to 
{Y} \otimes_F \mathfrak{F}$ for the base-extension
to $\mathfrak{F}$.
We call $X$ a
\emph{variety} over $F$ if $X$
is a
geometrically reduced algebraic scheme over $F$.

A variety $X$ over $\overline{K}$ is called a 
\emph{constant variety}
if there exists a variety $\widetilde{X}$ over $k$ with
$X = \widetilde{X} \otimes_k \overline{K}$.
Further, a subscheme $Y$ of $X$ is called a 
\emph{constant subscheme}
if
$Y = \widetilde{Y} \otimes_k \overline{K}$
for some subscheme $\widetilde{Y}$ of $\widetilde{X}$.
An abelian variety $A$ over $\overline{K}$
is called a \emph{constant abelian variety} if there exists an abelian variety
$\widetilde{A}$ over $k$ with $A = \widetilde{A} \otimes_k \overline{K}$
as abelian varieties.
Note that 
the group scheme structure of $A$
is required to be defined over $k$.

Let $A$ be an abelian variety over $\overline{K}$.
It is well known that
there exists a unique pair $\left( \widetilde{A}^{\overline{K}/k}, \Tr_A \right)$
consisting
of an abelian variety $\widetilde{A}^{\overline{K}/k}$ over $k$ and a homomorphism
$\Tr_A : \widetilde{A}^{\overline{K}/k} \otimes_{k} \overline{K} \to A$
of abelian varieties over $\overline{K}$
characterized by the property
that,
for any abelian variety $\widetilde{B}$ over $k$
and a homomorphism $\phi : \widetilde{B} \otimes_{k} \overline{K}
\to A$,
there exists a unique homomorphism $\Tr (\phi) :
\widetilde{B} \to 
\widetilde{A}^{\overline{K}/k} $ such that
$\phi$ factors as
$\phi = \Tr_A \circ \left( \Tr (\phi) \otimes_k \overline{K} \right)$.
This pair is called the \emph{$\overline{K}/k$-trace},
or simply the \emph{trace}, of $A$.
We sometimes call $\widetilde{A}^{\overline{K}/k}$ 
the $\overline{K}/k$-trace
by abuse of words.
See \cite{lang1} for more details.

\subsection{Geometric Bogomolov conjecture 
and known results} \label{sect:GBC}
We review
the geometric Bogomolov conjecture
for abelian varieties
and 
some
known results.
%
%
We begin by recalling the special subvarieties,
introduced in \cite{yamaki5}
and
used
to formulate the geometric Bogomolov conjecture for
abelian varieties.
They are
defined as
irreducible closed subvarieties
which are expressed as the sum
of the image of a constant closed subvariety in the trace
and a torsion subvariety.
To be precise, 
let $A$ be an abelian variety over $\overline{K}$
and
let $\left( \widetilde{A}^{\overline{K}/k} , \Tr_A \right)$
be the trace of $A$.
An
irreducible closed subvariety
$X$ of $A$ is called a \emph{special subvariety}
if
there exist a closed subvariety $\widetilde{Y}$ of 
$\widetilde{A}^{\overline{K}/k}$  and
a torsion subvariety $T \subset A$
such that 
$
X = T + \Tr_A {\left( \widetilde{Y} \otimes_k \overline{K} \right)}
$.%
\footnote{
In \cite{scanlon}, there is a similar but different notion
of special subvarieties due to Scanlon
(cf. \cite[Remark~7.3]{yamaki6}).}

In stating the 
conjecture,
it is convenient to
use the words of ``density
of small points'', which we are going to explain.
Let $L$ be an even ample line bundle on $A$.
Here $L$ is said to be even if $[-1]^{\ast} (L) =L$,
where $[-1] : A \to A$ is the homomorphism $a \mapsto -a$.
Let $\widehat{h}_{L}$ be the canonical height function over 
$A \left( \overline{K} \right)$
associated to $L$.
It is a semi-positive quadratic form 
on the group $A \left( \overline{K} \right)$.
We set, for $\epsilon > 0$,
\[
X ( \epsilon ; L) 
:=
\left\{
x \in X \left( \overline{K} \right)
\left|
\widehat{h}_{L} ( x) \leq \epsilon
\right.
\right\}
.
\]
It follows from \cite[Lemma~2.1]{yamaki5} that
whether or not
$X ( \epsilon ; L) $ is dense in $X$
for any $\epsilon > 0$
does not depend
on the choice of even ample $L$.
Therefore, it makes sense to say that
\emph{$X$ has dense small points}
if $X ( \epsilon ; L) $ is dense in $X$ for any $\epsilon > 0$
and for some 
(and hence any) even ample line bundle $L$ on $A$
(cf. \cite[Definition~2.2]{yamaki5}).

A point $x \in A \left( \overline{K} \right)$ is called a \emph{special point} if $\{ x \}$
is a special subvariety.
It is classically known that
a point is special if and only if it has height $0$
(cf. \cite[(2.5.4)]{yamaki5}).
This means that,
for an irreducible subvariety of dimension $0$,
being special is the same thing as having dense small points.
In the case of positive dimension also, 
it is verified
that
any special subvariety has dense small points
(cf. \cite[Corollary 2.8]{yamaki5}).
However,
it
is not known whether
the converse holds true or not in general.
The geometric Bogomolov conjecture for abelian varieties
is 
stating the converse:

\begin{Conjecture} [Geometric Bogomolov conjecture for
abelian varieties] \label{GBCforAV}
Let $A$ be an abelian variety over $\overline{K}$.
Then 
any irreducible closed subvariety
with dense small points
should be
a special subvariety.
\end{Conjecture}

Although Conjecture~\ref{GBCforAV} is not solved in full generality,
there are some results proved under some assumptions.
In \cite{gubler2},
Gubler proves that, if $A$ is totally degenerate 
at some place $v \in {M}_{\overline{K}}$,
then the conjecture holds true for $A$.
Here, $A$ is said to be totally degenerate at $v$
if the base-change of $A$ to 
$\overline{K}_v$
can be uniformized by
an algebraic torus
in the category of non-archimedean analytic spaces.
In \cite{yamaki5},
we generalize
this result for abelian varieties
allowing milder degenerations.
In \cite{yamaki6},
we moreover prove
that the conjecture holds 
for abelian varieties 
with ``nowhere degeneracy rank'' at most $1$.
(See
Theorem~\ref{TheoremD} below.)

Let us recall here the
notions of
maximal nowhere degenerate abelian subvariety of an abelian
variety and the nowhere degeneracy rank,
defined in \cite[Definition~7.10]{yamaki6}.
An abelian variety
$A$ over $\overline{K}$
is said to be \emph{non-degenerate at $v \in {M}_{\overline{K}}$}
if $A \otimes_{\overline{K}} \overline{K}_v$
is
the generic fiber of an abelian scheme over the ring of integers
of $\overline{K}_v$.
We say that $A$
is \emph{nowhere degenerate}
if
it is non-degenerate at any $v \in {M}_{\overline{K}}$.
We note
the fact that
there exists a unique \emph{maximal} abelian subvariety
$\mathfrak{m}$ of $A$
such that $\mathfrak{m}$ is nowhere degenerate,
where ``maximal'' means ``maximal  with respect to inclusion''.
This unique abelian subvariety is called the 
\emph{maximal nowhere degenerate abelian subvariety} of $A$.
Further, the \emph{nowhere-degeneracy rank of $A$}
is defined to be $\ndr (A) = \dim (\mathfrak{m})$.

With the use of nowhere degeneracy rank,
one of the main results of 
\cite{yamaki6} is stated as follows.

\begin{citeTheorem} [Theorem~D in \cite{yamaki6}] \label{TheoremD}
Let $A$ be an abelian variety over $\overline{K}$ with
$\ndr (A) \leq 1$.
Then,
any irreducible closed subvariety of $A$
with dense small points is a special subvariety.
\end{citeTheorem}

In \cite{yamaki6}, 
we
also pointed out that
the conjecture for an abelian variety
is equivalent to the conjecture for its
maximal nowhere degenerate abelian subvariety:

\begin{citeTheorem} [Theorem~E in \cite{yamaki6}]
\label{thm:reduction-GBCforMNDAS}
Let $A$ be an abelian variety over $\overline{K}$
with maximal
nowhere degenerate abelian subvariety
$\mathfrak{m}$.
Then 
the
geometric Bogomolov 
conjecture holds for $A$
if and only if
it
holds for $\mathfrak{m}$.
\end{citeTheorem}

Since
the geometric Bogomolov conjecture holds for elliptic curves,
Theorem~\ref{TheoremD} follows from 
Theorem~\ref{thm:reduction-GBCforMNDAS}.
By Theorem~\ref{thm:reduction-GBCforMNDAS}, 
the 
conjecture
for \emph{any} abelian variety
is 
reduced to 
that for
nowhere degenerate abelian varieties
(cf. \cite[Conjecture~7.22]{yamaki6}).


\subsection{Results and ideas} 
This paper includes two main results.
One 
is the following theorem.

\begin{Theorem} [Theorem~\ref{thm:main1}] \label{thm:main1intro}
Let $A$ be an abelian variety over $\overline{K}$
such that
$\dim \left( \widetilde{A}^{\overline{K}/k} \right) = \ndr (A)$.
Then the geometric Bogomolov conjecture holds for $A$.
\end{Theorem}

Theorem~\ref{thm:main1intro} generalizes Theorem~\ref{TheoremD}
(see Remark~\ref{rem:genthmD} how it generalize Theorem~\ref{TheoremD}).
We notice that the proof of Theorem~\ref{thm:main1intro} uses
Theorem~\ref{thm:reduction-GBCforMNDAS}.

The other theorem is the following,
where note that
$
\Image ( \Tr_A ) \subset \mathfrak{m}$
by \cite[Proposition~7.11]{yamaki6}.

\begin{Theorem} [Theorem~\ref{thm:main2}] \label{thm:main2intro}
Let $A$ be an abelian variety over $\overline{K}$
with
maximal nowhere degenerate abelian subvariety $\mathfrak{m}$
and let
$\mathfrak{t} := \Image (\Tr_A)$ be the image of the 
$\overline{K}/k$-trace homomorphism
of $A$.
Then the following statements are equivalent to each other.
\begin{enumerate}
\renewcommand{\labelenumi}{(\alph{enumi})}
\item
The geometric Bogomolov conjecture holds for $A$.
\item
The geometric Bogomolov conjecture holds for $\mathfrak{m}$.
\item
The geometric Bogomolov conjecture holds for
$\mathfrak{m} / \mathfrak{t}$.
\end{enumerate}
\end{Theorem}

This theorem includes Theorem~\ref{thm:reduction-GBCforMNDAS},
and
the new part is the equivalence between (c) and the others.
Since
(b) implies
(c) by
\cite[Lemma~7.7]{yamaki6},
the essential part is that (c) implies (b), in fact.
We remark that 
Theorem~\ref{thm:main2intro} leads us to
Theorem~\ref{thm:main1intro}.
(See Remark~\ref{rem:1from2}.)

Theorem~\ref{thm:main2intro}
is interesting because
it
shows that
Conjecture~\ref{GBCforAV} can be
reduced to the geometric Bogomolov
conjecture for a quite special class of abelian varieties,
that is,
for
nowhere degenerate abelian varieties with trivial trace
(cf. Remark~\ref{rem:reductiontonew1} and
Conjecture~\ref{conj:new1}.)

We briefly outline the idea of the proofs.
Let $A$ be an abelian variety over $\overline{K}$.
The starting point of the proofs 
is the fact proved 
in \cite
{gubler2}
that a closed subvariety of 
$A$ 
has
dense small points if and only if it has canonical height $0$
(cf. Proposition~\ref{prop:dense-height0}).
Using this fact
together with the description of the canonical height
in terms of intersection theory
over a model
(cf. Proposition~\ref{prop:model1}
and Lemma~\ref{lem:height-intersection}), 
we show 
Proposition~\ref{prop:GBCforconstant}.
This claims
that, when $A$ is a constant abelian variety,
a closed subvariety of canonical height $0$
is a constant subvariety and hence a special subvariety.
Then Theorem~\ref{thm:main1intro} follows from this proposition
and Theorem~\ref{thm:reduction-GBCforMNDAS}.

The proof of Theorem~\ref{thm:main2intro} needs more arguments,
which we are now going to explain.
Since
we may replace $A$ with an isogenous abelian variety
in proving the geometric Bogomolov conjecture
(cf. \cite[Corollary~7.6]{yamaki6}),
we find that
Theorem~\ref{thm:main2intro} is reduced to
Theorem~\ref{thm:height-torsion},
which claims
the following.
Let $A$ be a nowhere degenerate abelian variety over $\overline{K}$
with trivial trace
and let $B$ be a constant abelian variety.
Let
$X$ be a closed subvariety of $B \times A$
and let $Y$ and $T$ be the projections of $X$ to $B$ and $A$,
respectively.
Suppose that $X$ has dense small points
and
assume that the geometric Bogomolov conjecture holds for $A$.
Then, $Y$ is a constant subvariety, 
$T$ is a torsion subvariety, and $X = Y \times T$.

We give the outline of the proof of Theorem~\ref{thm:height-torsion}.
Let $A$, $B$, $X$, $Y$ and $T$ be as above.
Suppose that $X$ has dense small points.
Then
$Y$ has dense small points 
by \cite[Lemma~7.7]{yamaki6}
and hence has canonical height $0$ by Proposition~\ref{prop:dense-height0}.
Therefore by Proposition~\ref{prop:GBCforconstant}, 
$Y = \widetilde{Y} \otimes_k \overline{K}$
for some variety
$\widetilde{Y}$ 
over $k$.
It remains 
to show that $T$ is torsion and 
$X = Y \times T$.
To do that,
the ``relative height'' of $X \to Y$,
which will be introduced in \S~\ref{subsec:relativeheight}, 
plays a crucial role.
It is a function which assigns to each
point of a dense open subset of $\widetilde{Y}$
the canonical height of the fiber of $X \to Y$ over 
the corresponding point of $Y$.
In the proof of Theorem~\ref{thm:height-torsion},
we show in fact that,
if 
$X$ has dense small points,
then
the relative height function vanishes over a dense open subset of 
$\widetilde{Y}$.
Then
the 
geometric Bogomolov conjecture for $A$
implies
that
the geometric generic fiber of $X \to Y$ is a torsion subvariety,
and
arguments
using Chow's theorem (cf. Proposition~\ref{prop:Chow})
allow us to
conclude
that $T$ is a torsion subvariety
and $X = Y \times T$.

\subsection{Organization}
This article consists of six sections including this section
and an appendix.
In \S~\ref{sect:canonicalheight},
we recall canonical heights,
and 
describe them in terms of intersection theory over a model
when the abelian variety is nowhere degenerate.
Using the results there,
we show Theorem~\ref{thm:main1intro}
in \S~\ref{sect:GBCforConst}.
In \S~\ref{sect:family},
we investigate families 
of closed subvarieties 
of an abelian variety 
parameterized by a constant variety,
and we introduce
the relative heights.
In \S~\ref{sect:application},
applying the arguments in \S~\ref{sect:family}
to the setting of the geometric Bogomolov conjecture,
we
prove Theorem~\ref{thm:main2intro}.
In the appendix,
we prove a lemma 
concerning the trace of an abelian variety and the base-change,
which
is used in \S~\ref{sect:family}. 

\subsection*{Acknowledgments}
This research was done in part during 
my visit
 to 
the Institute of Mathematics of Jussieu
in September 2013,
which was
supported by ANR R\'egulateurs.
I thank Professor Vincent Maillot for his hospitality.
I thank Professor Walter Gubler for helpful comments
on an earlier draft of this paper.
I also thank the referee for valuable comments.
%
This work was partly supported by KAKENHI 21740012
and by KAKENHI 26800012.

\section{Canonical heights over nowhere degenerate
abelian varieties} \label{sect:canonicalheight}

The purpose of this section is to describe the canonical height
of closed subvarieties of nowhere degenerate abelian varieties
in terms of intersection theory on models.

\subsection{Canonical heights}

We briefly recall properties
of canonical heights
of closed subvarieties and cycles of abelian varieties.
We 
refer to \cite{gubler0, gubler2, gubler3} for more details.

Let $L$ be a line bundle on a variety $W$ over $\overline{K}$.
A \emph{metric
on 
$L $ at $v$}
means a collection of $\overline{K}_v$-norms $L (w) \to \RR$ for all $w \in 
W \left( \overline{K}_v \right)$,
where $L(w) := w^{\ast} (L)$ is the fiber of $L$ at $w$.
A \emph{metric on $L$} is a family $|| \cdot ||
= \{ || \cdot ||_{v} \}_{v \in {M}_{\overline{K}}}$
of metrics on $L$
at $v$ for all places $v \in M_{\overline{K}}$.
A line
bundle $L$ with a metric $|| \cdot ||$ is called a \emph{metrized line bundle},
denoted by $\overline{L} = \left( L ,  || \cdot || \right)$.

\begin{Example} [Algebraic metrics] \label{ex:algebraicmetric}
Let $K'$ be a finite extension of $K$
and let $\mathfrak{B}'$ be the normalization of 
$\mathfrak{B}$ in $K'$.
Let $f : \mathscr{W} \to \mathfrak{B}'$
be a proper morphism with geometric generic fiber $W$
and let $\mathscr{L}$ be a line bundle on $\mathscr{W}$ which
equals $L$ over $W$.
Then, it is known that an algebraic metric $|| \cdot ||_{\mathscr{L}}$ on $L$
is defined (cf. \cite[2.3]{gubler3}).
Here,
we do not recall what it is exactly
but
explain what it is like.
For any
$v \in M_{\overline{K}}$, 
let $R_v$ be the ring of integers of $\overline{K}_v$.
Let
$\Spec \left( \overline{K}_v \right)
\to
\mathfrak{B}'$
be the morphism
arising from the field extension $\overline{K}_v / K'$.
Since
$\mathfrak{B}'$ is proper over $k$,
this morphism
extends to a unique morphism $\Spec ( R_v ) \to \mathfrak{B}'$
by the valuative criterion.
Let $\mathscr{W}_{v} \to \Spec (R_{v})$
be the base-change of $f$ by this morphism.
Take any $w \in W \left( \overline{K}_v \right)$.
Since $f$ is proper,
there exists a unique section
$\sigma_w : \Spec (R_{v}) \to \mathscr{W}_{v}$
corresponding to $w$.
Note that
$\sigma_{w}^{\ast} ( 
\mathscr{L} \otimes_{\OO_{\mathfrak{B}'}}
R_v) $ a free $R_v$-module of rank $1$
and $\sigma_{w}^{\ast} ( 
\mathscr{L} \otimes_{\OO_{\mathfrak{B}'}}
R_v) \otimes_{R_v} \overline{K}_v = L (w)$.
To
a non-zero $s (w) \in L (w)$,
assign a non-negative number
\[
|| s(w)||_{\mathscr{L} , v}
:=
\inf
\left\{
|a^{-1}|_{v, \mathcal{H}} \in \RR_{>0}
\left| a \in \overline{K}_v^{\times}
,
\ 
a s(w) \in \sigma_{w}^{\ast} ( 
\mathscr{L} \otimes_{\OO_{\mathfrak{B}'}}
R_v) 
\right.
\right\}.
\]
This assignment 
defines a metric on 
$L$ at $v$
for each $v \in M_{\overline{K}}$
and hence a metric $|| \cdot ||_{\mathscr{L}} =
\left\{ || \cdot ||_{\mathscr{L} , v} \right\}_{v \in M_{\overline{K}}}$ on $L$.
This metric is
called the \emph{algebraic metric} associated to the model 
$(f , \mathscr{L})$.
\end{Example}

\begin{Example} [Canonical metrics] \label{ex:canonical}
Let $A$ be an abelian variety over $\overline{K}$.
For any $n \in \ZZ$, let $[n] : A \to A$ denote the endomorphism given by
$a \mapsto na$.
Let $L$ be a line bundle on $A$
and
assume that $L$ is even, that is,
$[-1]^{\ast} (L) \cong L$.
Let $n$ be an integer with $n \geq 2$.
It follows from the theorem of cube
that
there exists 
an isomorphism
$\phi : [n]^{\ast} (L) \to L^{\otimes n^{2}}$.
A metric $|| \cdot ||$
on $L$  
is called a \emph{canonical metric}
if 
$\phi$ induces an isometry
$[n]^{\ast} \left( \overline{L} \right)
\cong \overline{L}^{\otimes n^{2}}$,
where $\overline{L} := (L , || \cdot ||)$.
It is known that,
once the isomorphism $\phi$ is fixed, the canonical metric is determined
uniquely,
and that, for a different choice of isomorphisms,
the canonical metric changes only by
a non-zero constant multiple
(cf. \cite[Theorem~2.2]{zhang1-2} or \cite[Theorem~10.9]{gubler0}).
\end{Example}

To 
define the height with respect to metrized line bundles,
we need to 
focus on suitable metrics,
called \emph{admissible metrics},
studied in \cite{gubler0}.
We do not repeat the definition of them here
but
remark that any algebraic metric is known to be admissible.
Further, we also remark that canonical metrics over abelian varieties
are also admissible.
Thus we may take the heights
of closed subvarieties
 with respect to line bundles equipped
with these metrics.

Here is a remark concerning the compatibility
of metrics considered here with those developed in 
\cite{gubler0}.
Metrics here are those on
line bundles over an algebraic variety
and are considered only at the closed points of the variety.
On the other hand,
metrics
in \cite{gubler0}
are
those
on line bundles over the
analytic space associated to the algebraic
variety
and are considered at any point of the analytic space.
As far as working with admissible metrics, however,
we do not have to be serious about this difference.
In fact,
since
the admissible metrics are continuous metrics
on the analytic space and
the set of closed points is dense in the analytic space,
one can recover
all information of metrics on 
line bundles over the analytic space
from 
metrics over the closed points
of the given algebraic variety.

Let $\overline{L_0} , \ldots , \overline{L_d}$ be 
admissibly
metrized line bundles on a proper variety
$W$
and
let $Z$ be a $d$-dimensional cycle on a variety $W$.
Then the height
$
h_{\overline{L_0}, \ldots , \overline{L_d}} (Z)
$
of 
$Z$ 
with respect to 
$\overline{L_0} , \ldots , \overline{L_d}$
is defined in \cite[Definition~3.6]{gubler3}.
We do not recall
the definition of heights
because we do not need it
in the following arguments.
We
note that the heights that will be used in the arguments
are those
with respect to line bundles with algebraic metric
only,
which
can be described in terms of intersection theory
(cf. \cite[Theorem~3.5~(d)]{gubler3} or
Lemma~\ref{lem:height-intersection}).

Now, we consider an abelian variety $A$ over $\overline{K}$.
Let 
${L_0} , \ldots , {L_d}$
be even line bundles over $A$.
Fixing 
an isomorphism $[n]^{\ast} (L_i)
\to
L_i^{\otimes n^{2}}$ for each $i = 0 , \ldots , d$,
we obtain a canonically metrized line bundle
$\overline{L_i} = ( L_i , || \cdot ||_i)$
(cf. Example~\ref{ex:canonical}).
Then the \emph{canonical height} of a cycle $Z$ of dimension $d$
of $A$ with respect to
$L_0 , \ldots , L_d$ is defined to be
\[
\widehat{h}_{L_0 , \ldots , L_d} (Z)
:=
h_{\overline{L_0}, \ldots , \overline{L_d}} (Z).
\]
The canonical metrics 
on a line bundle
depend on the choice of isomorphisms 
$[n]^{\ast} (L_i)
\to
L_i^{\otimes n^{2}}$,
but it follows from the product formula that 
the canonical height does not.
Thus
the canonical height $\widehat{h}_{L_0 , \ldots , L_d} (Z)$
is well-defined.
If $L = L_0 = \cdots = L_d$
and there is no danger of confusion,
we simply write $\widehat{h}_{L} (Z)$
for
$\widehat{h}_{L_0 , \ldots , L_d} (Z)$.

For a closed subvariety $X$ of $A$ of pure dimension $d$,
let $[X]$ denote the corresponding cycle.
We write $\widehat{h}_{L_0 , \ldots , L_d} (X)
:=
\widehat{h}_{L_0 , \ldots , L_d} ([X])$,
called the
the canonical height of $X$
with respect to $L_0 , \ldots , L_d$.

The following proposition can be
found in  \cite{gubler2}.

\begin{citeProposition} [Corollary~4.4 in \cite{gubler2}] \label{prop:dense-height0}
Let $A$ be an abelian variety over $\overline{K}$,
let
$L$ be an even ample line bundle on $A$,
and
let $X$ be an irreducible closed subvariety
of $A$.
Then
$X$ 
has dense small points if and only if
$\widehat{h}_{L} (X) = 0$.
\end{citeProposition}

By Proposition~\ref{prop:dense-height0},
the geometric Bogomolov conjecture for $A$
is equivalent to the statement that
an irreducible closed subvariety $X$ of $A$
with
$\widehat{h}_{L} (X) = 0$ for some
even ample line bundle $L$ on $A$
is a special subvariety.

\subsection{Models of nowhere degenerate abelian varieties}

Let $W$ be a projective scheme over $\overline{K}$
and let $L$ be a line bundle on $W$.
Let $\mathfrak{U}$ be an open subset of $\mathfrak{B}'$.
A proper morphism $f : \mathscr{W} \to \mathfrak{U}$
with geometric generic fiber $W$
is called a \emph{model} of $W$ over $\mathfrak{U}$.
We note that,
in this terminology, \emph{it is not required that 
$W$ is dense in
$\mathscr{W}$}.
Let $\mathscr{L}$ be a line bundle on $\mathscr{W}$
such that the restriction of $\mathscr{L}$ to the
geometric generic fiber equals $L$.
Then the pair $(f, \mathscr{L})$ is called a \emph{model} of $(W,L)$
over $\mathfrak{U}$.

We construct a suitable model of nowhere degenerate abelian varieties.

\begin{Proposition} \label{prop:model1}
Let $A$ be a nowhere degenerate abelian variety over $\overline{K}$
and let $L$ be a line bundle on $A$.
Then, there exist a finite extension $K'$ of $K$,
a proper morphism $f : \mathscr{A} \to \mathfrak{B}'$,
where $\mathfrak{B}'$ is the 
normalization
of $\mathfrak{B}$ in $K'$,
and 
a line bundle $\mathscr{L}$ on $\mathscr{A}$
satisfying the following conditions.
\begin{enumerate}
\renewcommand{\labelenumi}{(\alph{enumi})}
\item
The pair $\left( f , \mathscr{L} \right)$
is a model of $(A ,L)$.
\item
There exists an open subset $\mathfrak{U} \subset \mathfrak{B}'$
with $\codim ( \mathfrak{B}' \setminus \mathfrak{U}
,  \mathfrak{B}' ) \geq 2$
such that
the restriction $f' : 
\mathscr{A}_{\mathfrak{U}} :=
f^{-1} (\mathfrak{U}) \to \mathfrak{U}$ 
of $f$ is 
an abelian scheme.
\item
Let $0_{f'}$ be the zero-section of the abelian
scheme $f'$ in (b).
Then $0_{f'}^{\ast} ( \mathscr{L} ) \cong \OO_{\mathfrak{U}}$.
\end{enumerate}
\end{Proposition}

\Proof
Let $K_0$ be a finite extension of $K$
such that $A$ and $L$ can be defined over $K_0$, that is,
$A = A_0 \otimes_{K_0} \overline{K}$ for some abelian variety 
$A_0$ over $K_0$ and $L$ is the base-change of a line
bundle on $A_0$.
Let $\mathfrak{B}_0$ be the normalization of $\mathfrak{B}$ in $K_0$.
Then 
there exist a dense open subset $\mathfrak{U}_0 \subset \mathfrak{B}_0$
and an abelian scheme $f_{0} : \mathscr{A}_0 \to \mathfrak{U}_0$
with zero-section $0_{f_{0}}$
having $A_0$ as its
generic fiber.

There exist only a finite number of points 
of $\mathfrak{B}_0$
of codimension $1$ which are 
not contained in $\mathfrak{U}_0$,
so that
let $v_1 , \ldots , v_m$ be all such points.
For each $i = 1 , \ldots , m$,
let $\OO_{i}$ be the stalk at $v_i$ of the structure sheaf of $\mathfrak{B}_0$.
It is a discrete valuation ring with fractional field $K_0$.
It follows from Grothendieck's semistable reduction theorem
(cf. \cite[Exp~IX, Th\'eor\`eme~3.6]{sga7})
that
there exists
a finite extension $K'$ of $K_0$
such that, for any $i = 1 , \ldots , m$, 
there exists a semistable model $\varphi_i$
of
$A_0 \otimes_{K_0} K'$ over
the integral closure $\OO_i'$ of $\OO_i$ in $K'$.
By the assumption of nowhere degeneracy of $A$,
these semistable models are abelian schemes.

Let $\nu : \mathfrak{B}' \to \mathfrak{B}_0$ 
be the normalization of $\mathfrak{B}_0$ in $K'$.
Note that 
$\Spec ( \OO_i' ) = \mathfrak{B}' \times_{\mathfrak{B}_0} \Spec ( \OO_i)$.
The abelian scheme
$\varphi_i
$
extends to an abelian scheme over a neighborhood of $\Spec (\OO_i')$
in $\mathfrak{B}'$.
Since
$\OO_i'$ is finite over $\OO_i$,
it follows that,
for each $i = 1, \ldots , m$,
there exist an open neighborhood $\mathfrak{U}_i$ of 
$v_i$ in $\mathfrak{B}_0$
and an abelian scheme $f_i' : \mathscr{A}_i ' \to
\nu^{-1} ( \mathfrak{U}_{i} )$
with zero-section $0_{f_i'}$
such that 
the generic fiber of 
$f_i' : \mathscr{A}_i ' \to
\nu^{-1} ( \mathfrak{U}_{i} )$
equals $A_{0} \otimes_{K_0} K'$.
On the other hand,
taking the base-change 
by $\nu^{-1} ( \mathfrak{U}_0 ) \to \mathfrak{U}_0$
of the abelian scheme $f_{0} : \mathscr{A}_0 
\to \mathfrak{U}_0$,
we obtain an abelian scheme
$f_0' : \mathscr{A}_0' 
 \to \nu^{-1} (\mathfrak{U}_0)$
with zero-section 
$0_{f'_{0}}$,
which has
$A_{0} \otimes_{K_0} K'$ as 
its 
generic fiber.
Thus we have
a family
\[
\Phi
:=
\left\{ \left( f_i' : \mathscr{A}_i '\to \nu^{-1} (\mathfrak{U}_i) ,
0_{f_{i}'}
\right) \right\}_{i = 0}^{m}
\]
of abelian schemes with generic fiber $A_{0} \otimes_{K_0} K'$.

It follows from the 
generalized
Weil extension lemma (cf. \cite[Proposition~1.3]{artin})
together with the valuative criterion of properness
that 
the isomorphism between the generic fibers of
the abelian schemes $\left( f_i' , 0_{f_{i}'}\right)$ and 
$\left( f_j' , 0_{f_{j}'} \right)$ 
extends to a unique isomorphism
between the restrictions
over $\nu^{-1} (\mathfrak{U}_i) \cap  \nu^{-1} (\mathfrak{U}_j)$
of the
abelian schemes,
or in other words,
the abelian schemes $\left( f_i' , 0_{f_{i}'}\right)$ and 
$\left( f_j' , 0_{f_{j}'} \right)$ 
coincide with each other over
$\nu^{-1} (\mathfrak{U}_i) \cap  \nu^{-1} (\mathfrak{U}_j)$.
Thus
the family 
$\Phi$
patches together to be 
an abelian scheme 
$f'_{+} : \mathscr{A}'_{+} \to \bigcup_{i = 0}^{m} \nu^{-1} (\mathfrak{U}_i)$
with geometric generic fiber $A = A_{0} \otimes_{K_0} \overline{K}$.

Since $\nu 
\left( \bigcup_{i = 0}^{m} \nu^{-1} (\mathfrak{U}_i)
\right) = \bigcup_{i = 0}^{m} \mathfrak{U}_i$ contains all
the points 
of $\mathfrak{B}_0$
of codimension $1$
and since $\nu$ is finite,
$\bigcup_{i = 0}^{m} \nu^{-1} (\mathfrak{U}_i)$ contains all the points 
of $\mathfrak{B}'$
of codimension $1$,
which means
$\codim \left( \mathfrak{B}' \setminus 
\bigcup_{i = 0}^{m} \nu^{-1} \left( \mathfrak{U}_i \right)
,
\mathfrak{B}' \right) \geq 2$.
Let $\mathfrak{U}$ be the regular locus of 
$\bigcup_{i = 0}^{m} \nu^{-1} (\mathfrak{U}_i)$.
Since $\mathfrak{B}'$ is regular in codimension $1$,
it follows that $\codim ( \mathfrak{B}' \setminus 
\mathfrak{U} , \mathfrak{B}') \geq 2$.

Let $f' : \mathscr{A}' \to \mathfrak{U}$ be the restriction
of $f'_{+}$ over $\mathfrak{U}$.
Then $f'$
is an abelian scheme with geometric generic fiber $A$.
Since $\mathfrak{U}$ is regular and $f'$ is smooth,
$\mathscr{A}'$ is regular.
It follows that
there exists a line bundle $\mathscr{L}'_{1}$ on $\mathscr{A}'$
such that the pair $(f' , \mathscr{L}'_1)$
is a model of $(A,L)$ over $\mathfrak{U}$.
Put $\mathscr{L}' :=
\mathscr{L}'_{1}
\otimes
(f')^{\ast} \left(0_{f'}^{\ast} \left( \mathscr{L}'_{1} 
\right)
\right)^{-1}$,
where $0_{f'}$ is the zero-section of $f'$.
Then we have a model $(f' , \mathscr{L}')$
such that
$0_{f'}^{\ast} \left( \mathscr{L}'
\right)$
is trivial.
Finally, by Nagata's embedding theorem
(see \cite[Theorem~5.7]{vojta}
for a scheme-theoretic version),
there exists a proper morphism $f : \mathscr{A} \to \mathfrak{B}'$
such that $\mathscr{A}'$ is an open dense subscheme of $\mathscr{A}$
and that $f|_{\mathscr{A}'} = f'$,
and an invertible sheaf $\mathscr{L}$ on $\mathscr{A}'$
such that $\mathscr{L} |_{\mathscr{A}'} = \mathscr{L}'$.
The pair $(f , \mathscr{L})$
satisfies the conditions (a) and (c),
and also satisfies (b) with the above $\mathfrak{U}$.
This proves the proposition.
\QED

\begin{Remark} \label{rem:modelofconstant}
Let $B$ be a constant abelian variety over $\overline{K}$,
and we take an abelian variety
$\widetilde{B}$ over $k$ with
$B = \widetilde{B} \otimes_k \overline{K}$.
Let $\widetilde{M}$ be an even ample line bundle on $\widetilde{B}$
and
set $M := \widetilde{M} \otimes_k \overline{K}$.
Let $\pr_{\mathfrak{B}} : \widetilde{B} \times_{\Spec (k)} \mathfrak{B} \to \mathfrak{B}$
and $q : 
\widetilde{B} \times_{\Spec (k)} \mathfrak{B} \to \widetilde{B}$
be the canonical projections.
Then the pair 
$\left( \pr_{\mathfrak{B}} , q^{\ast} \left( \widetilde{M} \right) \right)$
is a model of $(B , M)$ satisfying the conditions 
in Proposition~\ref{prop:model1}.
\end{Remark}

\subsection{Height and intersection on a model}

In this subsection,
we 
describe 
the 
canonical heights of pure dimensional closed subschemes
of a nowhere degenerate abelian variety
in terms of intersection theory
over models.

Let $A$ be an abelian variety over $\overline{K}$
with an even ample line bundle
$L$.
Assume that $A$ is nowhere degenerate,
and 
we take a model $(f : \mathscr{A} \to \mathfrak{B}' , \mathscr{L})$
of $(A , L)$
satisfying the conditions in
Proposition~\ref{prop:model1},
where
$\mathfrak{B}'$ is the normalization of $\mathfrak{B}$ in 
some finite extension $K'$ of $K$.
Let
$\mathfrak{U}$ be an open subset of $\mathfrak{B}'$
as in (b) of Proposition~\ref{prop:model1}.
Then the restriction $f' : \mathscr{A}_{\mathfrak{U}} \to \mathfrak{U}$,
where $\mathscr{A}_{\mathfrak{U}} := f^{-1} (\mathfrak{U})$,
is an abelian scheme with zero-section $0_{f'}$ by (c).

It follows from (c) of Proposition~\ref{prop:model1} that,
for any $n \in \NN$,
there exists an isomorphism
\addtocounter{Claim}{1}
\begin{align} \label{eq:model-linebundle}
[n]^{\ast} \left( \mathscr{L}|_{\mathscr{A}_{\mathfrak{U}}}
\right) \to \left( \mathscr{L}|_{\mathscr{A}_{\mathfrak{U}}}
\right)^{\otimes n^{2}}
,
\end{align}
where $[n] : \mathscr{A}_{\mathfrak{U}} \to \mathscr{A}_{\mathfrak{U}}$
is the $n$-times endomorphism.
Indeed, 
since $L$ is even,
there exists an isomorphism 
$
[n]^{\ast} \left( L \right) \to 
L^{\otimes n^{2}}
$, which extends to an isomorphism
$[n]^{\ast} \left( \mathscr{L} |_{\mathfrak{U}}
\right) \to \left( \mathscr{L} |_{\mathfrak{U}}
\right)^{\otimes n^{2}} \otimes f^{\ast} ( \mathscr{M} )$
for some line bundle $\mathscr{M}$ on $\mathfrak{U}$.
Then we have
an isomorphism
\[
0_{f'}^{\ast}
\left(
[n]^{\ast} \left( \mathscr{L} |_{\mathfrak{U}}
\right)
\right) \to 
0_{f'}^{\ast}
\left(
\left( \mathscr{L} |_{\mathfrak{U}}
\right)^{\otimes n^{2}} \otimes f^{\ast} ( \mathscr{M} )
\right)
=
0_{f'}^{\ast} \left( \mathscr{L} |_{\mathfrak{U}} \right)^{\otimes n^{2}}
\otimes \mathscr{M}
.
\]
By the condition (c),
we have $0_{f'}^{\ast}
\left(
[n]^{\ast} \left( \mathscr{L} |_{\mathfrak{U}}
\right)
\right) = [n]^{\ast} \left( 0_{f'}^{\ast}
\left( \mathscr{L} |_{\mathfrak{U}} \right) \right) \cong \OO_{\mathfrak{U}}$
and $0_{f'}^{\ast} \left( \mathscr{L} |_{\mathfrak{U}}
\right)^{\otimes n^{2}} \cong \OO_{\mathfrak{U}}$,
which lead us to
$\mathscr{M} \cong \OO_{\mathfrak{U}}
$.
This shows the existence of an isomorphism (\ref{eq:model-linebundle}).

Let $|| \cdot ||_{\mathscr{L}}$ be the algebraic metric on $L$
associated to the model $(f , \mathscr{L})$
(cf. Example~\ref{ex:algebraicmetric}).
Noting
$\codim ( \mathfrak{B}' \setminus \mathfrak{U}
,\mathfrak{B}') \geq 2$,
we see 
that
the isomorphism (\ref{eq:model-linebundle})
gives rise to an isometry 
\[
[n]^{\ast} (L, || \cdot ||_{\mathscr{L},v})
\cong
(L, || \cdot ||_{\mathscr{L},v})^{\otimes n^{2}}
\]
for each $v \in M_{\overline{K}}$
(cf. Example~\ref{ex:algebraicmetric}).
Thus
the 
algebraic 
metric $|| \cdot ||_{\mathscr{L}}$
is a canonical metric.

The heights of cycles with respect to algebraically metrized
line bundles are described in terms of intersection theory,
and
so are
the canonical heights.
To be precise,
let $X$ be a closed subscheme of $A$ of pure dimension $d$.
Replacing $K'$ with a 
finite extension if necessary,
we assume that $X$ can be defined over $K'$.
Let $\mathscr{X}$ be the closure of $X$ in $\mathscr{A}$.
Note that $X$ is the geometric generic fiber of $\mathscr{X}$.
Let $[\mathscr{X}]$ denote the cycle corresponding to $\mathscr{X}$.
Let $\mathcal{H}'$ be the pull-back of $\mathcal{H}$
to $\mathfrak{B}'$.
Since $|| \cdot ||_{\mathscr{L}}$ is a canonical metric,
we have, by \cite[Theorem~3.5~(d)]{gubler3},
\addtocounter{Claim}{1}
\begin{align} \label{eq:canonicalheight}
\widehat{h}_{L} (X)
=
\frac{
\deg_{\mathcal{H}'}
f_{\ast}
\left(
\cherncl_1 ( \mathscr{L} )^{\cdot (d+1)} \cdot [ \mathscr{X} ]
\right)
}{[K' : K]}
,
\end{align}
where 
$\deg_{\mathcal{H}'}$ means the degree 
of the $(b-1)$-dimensional cycle on $\mathfrak{B}'$
with respect to $\mathcal{H}'$
and $[K' : K]$ is the extension degree.

In 
equality
(\ref{eq:canonicalheight}),
the model $f |_{\mathscr{X}} : \mathscr{X} \to \mathfrak{B}'$
is assumed to be the \emph{closure} of $X$,
but
this is 
too strong to use in the latter arguments.
In fact,
we verify
equality 
(\ref{eq:canonicalheight})
under a milder assumption on $\mathscr{X}$
as follows.

\begin{Lemma} \label{lem:height-intersection}
Let $A$, $L$, $X$, $f : \mathscr{A} \to \mathfrak{B}'$,
$\mathscr{L}$, and $\mathcal{H}'$ be as above.
Let $\mathscr{X}'$ be a 
pure dimensional
closed subscheme of $\mathscr{A}$
such that the restriction $f|_{\mathscr{X}'} : \mathscr{X}' \to \mathfrak{B}'$
has geometric generic fiber $X$
and
is flat over any point of codimension $1$ of $\mathfrak{B}'$.
Then we have
\[
\widehat{h}_{L} (X)
=
\frac{
\deg_{\mathcal{H}'}
f_{\ast}
\left(
\cherncl_1 ( \mathscr{L} )^{\cdot (d+1)} \cdot [\mathscr{X}']
\right)}{[K' : K]}
,
\]
where 
$\deg_{\mathcal{H}'}$ means the degree 
of the $(b-1)$-dimensional cycle on $\mathfrak{B}'$
with respect to $\mathcal{H}'$.
\end{Lemma}

\Proof
Let $\mathscr{X}$ be the closure of $X$ in $\mathscr{A}$.
Then there exists an effective
cycle $[\mathscr{V}]$ on $\mathscr{A}$ of dimension $d+b$
such that
$[\mathscr{X}'] = [\mathscr{X}] + [\mathscr{V}]$.
Since $\mathscr{X}'$ is flat over any point of $\mathfrak{B}'$
of codimension $1$, $\mathscr{X}'$ coincides with $\mathscr{X}$
over any point of $\mathfrak{B}'$
of codimension $1$.
Therefore, any irreducible component of
$f ( \mathscr{V} )$ has dimension not greater than
$
b-2$.
Since $\cherncl_1 ( \mathscr{L} )^{\cdot (d+1)} \cdot [\mathscr{V}]$
is a cycle class of dimension $b-1$
and can be represented by a cycle
supported in $\mathscr{V}$,
it follows that $f_{\ast}
\left(
\cherncl_1 ( \mathscr{L} )^{\cdot (d+1)} \cdot [\mathscr{V}]
\right) = 0$.
This concludes
$f_{\ast}
\left(
\cherncl_1 ( \mathscr{L} )^{\cdot (d+1)} \cdot [\mathscr{X}']
\right)
=
f_{\ast}
\left(
\cherncl_1 ( \mathscr{L} )^{\cdot (d+1)} \cdot [\mathscr{X}]
\right)$,
and thus,
by (\ref{eq:canonicalheight}),
the desired equality is obtained.
\QED

\section{Constant abelian varieties and the geometric Bogomolov conjecture}
\label{sect:GBCforConst}

\subsection{Proof of Theorem~\ref{thm:main1intro}}
In this subsection, we show 
Proposition~\ref{prop:GBCforconstant},
which indicates
that Conjecture~\ref{GBCforAV}
holds for abelian varieties which are isogenous to constant
abelian varieties
(cf. Remark~\ref{rem:GBCconstant}).
Using this proposition,
we prove
Theorem~\ref{thm:main1intro},
which states
that the geometric Bogomolov conjecture
holds for any abelian variety whose nowhere degeneracy rank equals
the dimension of the trace.

We begin with a preliminary argument.
We will
use the following lemma on intersection theory.

\begin{Lemma} \label{lem:representingcycle}
Let $\mathscr{X}$ be an irreducible proper variety over $k$, 
let $\mathscr{Y}$ be 
a closed subscheme of $\mathscr{X}$ 
of pure dimension $d$,
and
let $\left[ \mathscr{Y} \right]$ 
denote the corresponding cycle.
Let $\mathscr{L}$ be a line bundle and let 
$\mathscr{M}$ be a
sublinear system of the complete linear system
associated to $\mathscr{L}$.
Suppose that $\mathscr{M}$ is base-point free.
Let $e$ be a positive integer.
Then, for general $\mathscr{D}_1 , \ldots , \mathscr{D}_e \in 
\mathscr{M}$,
the cycle $[\mathscr{D}_1 \cap \cdots \cap \mathscr{D}_e \cap 
\mathscr{Y}]$
represents the cycle class
$ \cherncl_1 ( \mathscr{L} )^{\cdot e} \cdot [\mathscr{Y}]$.
\end{Lemma}

\Proof
Since $\mathscr{M}$ is base-point free,
a general
$\mathscr{D}_1 \in \mathscr{M}$ 
is away from any associated point of $\mathscr{Y}$.
Then,
at any point of the support of $\mathscr{D}_1$,
a local equation of $\mathscr{D}_1$ is $\OO_{\mathscr{Y}}$-regular.
Since $\mathscr{M}$ is base-point free again,
a general $\mathscr{D}_{2} \in \mathscr{M}$ 
is away from 
any associated point of $\mathscr{D}_1 \cap \mathscr{Y}$.
Then a local equation of 
$\mathscr{D}_{2}$ is $\OO_{\mathscr{D}_1 \cap \mathscr{Y}}$-regular.
Repeating this process,
we obtain, by induction, 
general $\mathscr{D}_1 , \ldots , \mathscr{D}_e$ such that
the sequence of their local equations is $\OO_{\mathscr{Y}}$-regular.
By \cite[Example~2.4.8]{fulton},
therefore,
the cycle of $\mathscr{D}_1 \cap \cdots \cap \mathscr{D}_e \cap 
\mathscr{Y}$ represents
the cycle class $
\cherncl_1 ( \mathscr{L} )^{\cdot e} \cdot [\mathscr{Y}]$.
\QED

Now, we show the key proposition
to the proof of Theorem~\ref{thm:main1intro}.

\begin{Proposition} \label{prop:GBCforconstant}
Let $\widetilde{B}$ be an abelian variety over $k$
and let $\widetilde{M}$ be an even ample line bundle on $\widetilde{B}$.
Set
$B := \widetilde{B} \otimes_k \overline{K}$
and $M := \widetilde{M} \otimes_k \overline{K}$.
Let
$Y$ be an irreducible closed subvariety of $B$ of dimension $d$.
Suppose that
$\widehat{h}_{M} (Y) = 0$.
Then $Y$ is a constant subvariety.
\end{Proposition}

\Proof
Since the canonical height is multilinear
on line bundles to which it is associated
(cf. \cite[Theorem~11.18~(a)]{gubler0}),
we have
$\widehat{h}_{M^{\otimes n}} (Y) = n^{d+1} \widehat{h}_{M} (Y)$
for any $n \in \NN$,
so that
we may and do assume that 
$\widetilde{M}$ and hence $M$ are 
very ample.

Let $K'$ be a finite extension of $K$
such that $Y$ can be defined over $K'$
and let $\mathfrak{B}' \to \mathfrak{B}$ be the normalization of $\mathfrak{B}$
in $K'$.
We set $\mathscr{B} := \widetilde{B} \times_{\Spec (k)} \mathfrak{B}'$
and let $f : \mathscr{B} \to \mathfrak{B}'$ 
denote the second projection,
which is an abelian scheme.
Further, let $\mathscr{M}$ be the pull-back of $\widetilde{M}$
by the canonical projection $\pr_{\widetilde{B}} : 
\mathscr{B} \to \widetilde{B}$.
Then the pair $( f , \mathscr{M} )$
is a model of $(B , M)$ such that the conditions in Proposition~\ref{prop:model1} are satisfied
(cf. Remark~\ref{rem:modelofconstant}).
Let $\mathscr{Y}$ be the closure of $Y$ in $\mathscr{B}$.
Then the geometric generic fiber of
$f |_{\mathscr{Y}} : \mathscr{Y} \to \mathfrak{B}'$
equals $Y$.

Set $\widetilde{Y} := \pr_{\widetilde{B}} ( \mathscr{Y} )$.
We are going to show that $Y = \widetilde{Y} \otimes_k \overline{K}$.
Let $\mathfrak{U} \subset \mathfrak{B}'$ be a dense open subset
such that $f|_{\mathscr{Y}}$ is flat over $\mathfrak{U}$.
We put $\mathscr{V} := f^{-1} (\mathfrak{U}) \cap \mathscr{Y}$.
Note that
$f |_{\mathscr{V}} : \mathscr{V} \to \mathfrak{U}$ is a proper flat
morphism of relative dimension $d$.
Let $\left| \widetilde{M} \right|$
be the complete
linear system 
associated to $\widetilde{M}$.
Since $\widetilde{M}$ is very ample,
it follows that,
for
general members $\widetilde{D}_{1} , \ldots , \widetilde{D}_{d+1}
\in \left| \widetilde{M} \right|$,
if we write
$\mathscr{D}_i := \pr_{\widetilde{B}}^{-1} \left( 
\widetilde{D}_i
\right)$ for $i = 1 ,\ldots , d+1$,
then
the restriction
$
\mathscr{D}_1 \cap \cdots \cap \mathscr{D}_{d} 
\cap \mathscr{V}
\to \mathfrak{U}
$
of $f$
is finite and surjective,
and
furthermore,
by Lemma~\ref{lem:representingcycle},
the cycles of the closed subschemes
$
\mathscr{D}_1 \cap \cdots \cap \mathscr{D}_{d} 
\cap
\mathscr{Y}
$
and
$
\mathscr{D}_1 \cap \cdots \cap \mathscr{D}_{d+1} 
\cap
\mathscr{Y}
$
represent the cycle classes
$
\cherncl_1 \left( \mathscr{M} \right)^{\cdot d}
\cdot
[\mathscr{Y}]
$
and
$
\cherncl_1 \left( \mathscr{M} \right)^{\cdot (d + 1)}
\cdot
[\mathscr{Y}]
$
respectively.

For such $\mathscr{D}_1 , \ldots , \mathscr{D}_{d+1}$,
we have
\addtocounter{Claim}{1}
\begin{align} \label{eq:emptyintersection}
\mathscr{D}_1 \cap \cdots \cap \mathscr{D}_{d+1} \cap \mathscr{V}
= \emptyset.
\end{align}
Indeed, 
let $\mathscr{Z}$
be 
the closure of $\mathscr{D}_1 \cap \cdots \cap \mathscr{D}_{d+1} \cap \mathscr{V}$.
Then $[\mathscr{Z}]$ 
and $[\mathscr{D}_1 \cap \cdots \cap \mathscr{D}_{d+1} \cap \mathscr{Y}]
-
[\mathscr{Z}]$
are
effective cycles\footnote{Possibly the zero-cycle.} 
of dimension $b-1$,
where we remark $b = \dim ( \mathfrak{B}' )$.
Therefore we have
\begin{align*}
\deg_{\mathcal{H}'}
f_{\ast}
\left(
\cherncl_1 \left( \mathscr{M} \right)^{\cdot (d + 1)}
\cdot
[\mathscr{Y}]
\right)
=
\deg_{\mathcal{H}'}
f_{\ast}
\left(
\mathscr{D}_1 \cap \cdots \cap \mathscr{D}_{d+1} \cap \mathscr{Y}
\right)
\geq
\deg_{\mathcal{H}'}
f_{\ast}
[\mathscr{Z}]
\geq
0
,
\end{align*}
where 
$\mathcal{H}'$ is the pull-back of $\mathcal{H}$ by 
the finite morphism
$\mathfrak{B}' \to \mathfrak{B}$,
and
$\deg_{\mathcal{H}'}$ means the degree 
of the $(b-1)$-dimensional cycle on $\mathfrak{B}'$
with respect to $\mathcal{H}'$.
On the other hand,
by (\ref{eq:canonicalheight}) or Lemma~\ref{lem:height-intersection},
we have
$
\deg_{\mathcal{H}'}
f_{\ast}
\left(
\cherncl_1 \left( \mathscr{M} \right)^{\cdot (d + 1)}
\cdot
[\mathscr{Y}]
\right) = 
[K' : K]
\widehat{h}_{M} (Y)
$,
and hence
\[
\deg_{\mathcal{H}'}
f_{\ast}
\left(
\cherncl_1 \left( \mathscr{M} \right)^{\cdot (d + 1)}
\cdot
[\mathscr{Y}]
\right) = 0
\]
by assumption.
It follows that
$\deg_{\mathcal{H}'}
f_{\ast}
[\mathscr{Z}]
=0$.
Since $f_{\ast}
[
\mathscr{Z}]
$
is an effective cycle,
this implies that
the support $\Supp 
\left(
f_{\ast}
[
\mathscr{Z}
]
\right)
$
of $f_{\ast}
[
\mathscr{Z}
]$
is empty.
Since
any irreducible component of $\Supp ([\mathscr{Z}])$
is generically finite over $\mathfrak{B}'$,
we have $\Supp \left(
f_{\ast}
[
\mathscr{Z}
]
\right) =f \left( \Supp ([\mathscr{Z}]) \right)$.
Thus
$\Supp ([\mathscr{Z}]) = \emptyset$ follows,
which concludes
$
\mathscr{D}_1 \cap \cdots \cap \mathscr{D}_{d+1} \cap \mathscr{V}
= \emptyset$.

Next, we claim that $\dim \left( \widetilde{Y} \right)
= d$.
Since $\dim \left( \mathscr{Y} \right) = d + b$
and the generic fiber $Y$ of $f|_{\mathscr{Y}} 
: \mathscr{Y} \to \mathfrak{B}'$ has dimension $d$,
the dimension counting
shows 
$
\dim \left( \widetilde{Y} \right) =
\dim ( \pr_{\widetilde{B}} ( \mathscr{Y} ) )
\geq d$,
so that it suffices to show $\dim \left( \widetilde{Y} \right) \leq d$.
By (\ref{eq:emptyintersection}),
we have
\addtocounter{Claim}{1}
\begin{align} \label{eq:emptyintersection2}
\pr_{\widetilde{B}}
\left(
\mathscr{D}_1 \cap \cdots \cap \mathscr{D}_{d+1} \cap \mathscr{V}
\right)
=
\emptyset
\end{align}
for general  $\widetilde{D}_{1} , \ldots , \widetilde{D}_{d+1}
\in \left| \widetilde{M} \right|$,
where we recall
$\mathscr{D}_i := \pr_{\widetilde{B}}^{-1} \left( 
\widetilde{D}_i
\right)$ for $i = 1 ,\ldots , d+1$.
Since $\mathscr{V}$ is a dense open subset of $\mathscr{Y}$
and $\pr_{\widetilde{B}} |_{\mathscr{Y}} : \mathscr{Y} \to 
\widetilde{B}$ is surjective,
there
exists 
a dense open subset
$\widetilde{W}$ of 
$\widetilde{Y}$ with
$\widetilde{W} \subset  \pr_{\widetilde{B}} ( \mathscr{V})$.
Then (\ref{eq:emptyintersection2})
implies
$\widetilde{D}_1 \cap \cdots \cap \widetilde{D}_{d+1}
\cap \widetilde{W} = \emptyset$
for general  $\widetilde{D}_{1} , \ldots , \widetilde{D}_{d+1}
\in \left| \widetilde{M} \right|$.
Since $\widetilde{M}$ is very ample and
$\widetilde{W}$ is a dense open subset of $
\widetilde{Y}
$,
it follows that $
\dim \left( \widetilde{Y} \right) \leq d
$.
Thus, we find $\dim \left( \widetilde{Y} \right) = d$.

Since $\widetilde{Y} = \pr_{\widetilde{B}} ( \mathscr{Y} )$,
we
have $\mathscr{Y} \subset \widetilde{Y} \times_{\Spec (k)} \mathfrak{B} 
$.
Since both $\mathscr{Y}$
and
$\widetilde{Y} \times_{\Spec (k)} \mathfrak{B}$
are irreducible varieties of dimension $d+b$,
we conclude that
$\mathscr{Y} = \widetilde{Y} \times_{\Spec (k)} \mathfrak{B} $.
This proves
that $Y = \widetilde{Y} \otimes_k \overline{K}$,
and thus
$Y$ is a constant subvariety of $B$.
\QED

\begin{Remark} \label{rem:GBCconstant}
It follows from the definition of special subvarieties
that
any constant closed subvariety of 
a constant abelian variety is a special subvariety.
Taking into account of
Proposition~\ref{prop:dense-height0},
we then find that
Proposition~\ref{prop:GBCforconstant}
verifies the geometric Bogomolov conjecture
for constant abelian varieties.
\end{Remark}

Recall that,
for an abelian variety $A$
over $\overline{K}$, 
$\left( \widetilde{A}^{\overline{K}/k} , \Tr_A \right)$ denotes the 
$\overline{K}/k$-trace of $A$.

\begin{Remark} \label{rem:trace-ndr}
The trace homomorphism 
$\Tr_A$ is a finite morphism
(cf. \cite[Lemma~1.4]{yamaki5}).
It is known that $\Tr_A$ factors through $\mathfrak{m}$ by 
\cite[Proposition~7.11]{yamaki6},
where $\mathfrak{m}$ is the 
maximal nowhere degenerate abelian subvariety of $A$.
Thus
$\dim \left( \widetilde{A}^{\overline{K}/k} \right) \leq \ndr (A)$
holds.
\end{Remark}

Now we establish Theorem~\ref{thm:main1intro}
as a consequence of Proposition~\ref{prop:GBCforconstant}.

\begin{Theorem} [Theorem~\ref{thm:main1intro}] \label{thm:main1}
Let $A$ be an abelian variety over $\overline{K}$.
Assume
$\dim \left( \widetilde{A}^{\overline{K}/k} \right) = \ndr (A)$.
Then the geometric Bogomolov conjecture holds for $A$.
\end{Theorem}

\begin{Pf}
Let $\mathfrak{m}$ be the maximal nowhere degenerate abelian subvariety
of $A$.
We see from
Remark~\ref{rem:trace-ndr}
that
$\Tr_A$ factors as a finite homomorphism $\widetilde{A}^{\overline{K}/k} 
\otimes_k \overline{K} \to \mathfrak{m}$,
and
it follows from the assumption that this
homomorphism
is an isogeny.
By Remark~\ref{rem:GBCconstant},
which is a re-interpretation of Proposition~\ref{prop:GBCforconstant},
the geometric Bogomolov conjecture holds for $\widetilde{A}^{\overline{K}/k} 
\otimes_k \overline{K}$,
so that,
by \cite[Corollary~7.6]{yamaki6},
the conjecture holds also for $\mathfrak{m}$.
By Theorem~\ref{thm:reduction-GBCforMNDAS},
the conjecture holds for $A$.
\end{Pf}

\begin{Remark} \label{rem:genthmD}
%
%
%
Theorem~\ref{thm:main1intro},
which is Theorem~\ref{thm:main1} also, 
generalizes Theorem~\ref{TheoremD}.
Indeed, suppose that $\ndr (A) \leq 1$.
If $\ndr (A) = 0$, then Remark~\ref{rem:trace-ndr}
tells us
$\dim \left( \widetilde{A}^{\overline{K}/k} \right) = \ndr (A) = 0$,
so that Theorem~\ref{TheoremD} follows from Theorem~\ref{thm:main1}
trivially.
Suppose that $\ndr (A) = 1$.
This means that the maximal nowhere degenerate abelian
subvariety $\mathfrak{m}$
of $A$ has dimension $1$.
We remark 
here 
that
any nowhere degenerate abelian variety of dimension
$1$ is a constant abelian variety,
which follows from the well-known fact that the moduli
space of elliptic curves is an affine line.
Therefore,
$\mathfrak{m}$ is a constant abelian variety,
and it follows from
the universality of the trace
that
$\dim \left( \widetilde{A}^{\overline{K}/k} \right) \geq \dim (\mathfrak{m})$.
By Remark~\ref{rem:trace-ndr},
we then conclude that
$\dim \left( \widetilde{A}^{\overline{K}/k} \right)
= \dim ( \mathfrak{m} ) = \ndr (A)$.
Thus
Theorem~\ref{TheoremD} follows from Theorem~\ref{thm:main1}.
\end{Remark}

\subsection{Special subvarieties of constant abelian varieties}
First, we note the
following proposition, which
is a slight generalization of Chow's theorem
\cite[II,\S~1~Theorem~5]{lang1}
and will be significantly used
later 
to prove Theorem~\ref{thm:main2intro}.

\begin{Proposition} \label{prop:Chow}
Let $F$ be an algebraically closed field
and let
$\mathfrak{F} /F$ be a field extension
with $\mathfrak{F}$ algebraically closed.
Let $A$ be an abelian variety over $F$
and
let
$Z$ be a torsion subvariety of $
A \otimes_{F} \mathfrak{F}$.
Then there exists a torsion subvariety $T$ of $A$
such that $Z = T \otimes_{F} \mathfrak{F}
$.
\end{Proposition}

\Proof
Since $F$ is algebraically closed,
the torsion points of $\left( A \otimes_{F} \mathfrak{F} \right)
\left( \mathfrak{F} \right)$ coincide with
those of $A \left( F \right)$.
Then the assertion follows from
Chow's theorem \cite[II,\S~1~Theorem~5]{lang1}.
\QED

Proposition~\ref{prop:Chow}
is used in the following remark.
This remark is not necessary in the sequel,
but it
is worth mentioning
because
it tells us a basic fact which
characterizes
special subvarieties
of constant abelian varieties.

\begin{Remark}
Let $B
$ 
be a constant abelian variety over $\overline{K}$
and 
let $Y$ be an irreducible closed subvariety of $B$.
We then note that $Y$ is a special subvariety
if and only if $Y$ is a constant subvariety.
Indeed,
the ``if'' part
is noted in Remark~\ref{rem:GBCconstant}.
As for the ``only if'' part,
we have two different proofs:
one is the proof
using
Proposition~\ref{prop:GBCforconstant}
together with some basic facts on special subvarieties;
the other is the proof
just
using the definition of special subvarieties and
Proposition~\ref{prop:Chow}. 
The details are left to the reader.
\end{Remark}

\section{Family of closed subvarieties over a constant variety}
\label{sect:family}

In this section, we investigate
a family of closed subvarieties of an abelian variety
which is
parameterized by a constant variety.

\subsection{Abelian subschemes and their translates}

Let $p : \mathcal{A} \to U$ be an abelian scheme over
a noetherian scheme $U$ with
addition
$\alpha : \mathcal{A} \times_U \mathcal{A} \to \mathcal{A}$
and
zero-section $0_{p}$.
Let $X$ be a 
closed 
subscheme of $\mathcal{A}$ and let $p|_{X} : X \to U$
be the restriction of $p$.
We call $p|_{X} : X \to U$
a \emph{closed subgroup scheme} of $p$
if
the group structure of $p$
restricts to a group structure on $p|_{X}$,
that is,
we have
$\alpha \left( X \times_U X \right) \subset X$,
$0_{p} (U) \subset X$, and
$-X = X$,
where $-X$ is the image of $X$ by the ``minus'' morphism of the
abelian scheme $p$.

The following lemma shows that
the condition $-X = X$
can be omitted for $X$ to be a subgroup scheme.


\begin{Lemma} \label{lem:abeliansubscheme}
Let $
p : \mathcal{A} \to U$
be an abelian scheme
with addition 
$\alpha : \mathcal{A} \times_U \mathcal{A} \to \mathcal{A}$
and zero-section
$0_{p}$.
Let 
$X$ be a 
closed 
subscheme of $\mathcal{A}$
such that 
$\alpha \left( X \times_U X \right) \subset X$
and
$0_{p} (U) \subset X$.
Then,
$p|_{X} : X \to U$ is 
a closed subgroup scheme of $p$.
\end{Lemma}

\Proof
It suffices to show that
$-X = X$.
Let $S$ be any $U$-scheme and let
$s \in X(S)$ be any $S$-valued point of $X$.
Since the zero $0_{p,S}$ of the group $\mathcal{A} (S)$
sits in $X(S)$ by assumption,
the addition $
X (S) \times X(S) \to
X (S)$ is surjective.
Therefore,
there exists an
$s' \in X(S)$ such that
$s+s' = \alpha ( s, s' ) = 0_{p,S}$,
which shows that
$-s = s' \in X(S)$.
Thus we have
$- X(S) \subset X(S)$.
Since $S$ is an arbitrary $U$-scheme, 
this means that
$-X \subset X$.
We then have
$X = - (-X) \subset - X$,
and thus $-X = X$ as required.
\QED

Let
$p |_X : X \to U$ be a closed subgroup scheme of an abelian scheme
$p : \mathcal{A} \to U$.
We call 
$p|_{X}$ 
an \emph{abelian subscheme} of $p$
if it
is a smooth morphism with geometrically connected fibers.
Then $p|_X$ itself is an abelian scheme 
in a natural way.

The following lemma shows that, under a certain condition,
if there exists a
geometric fiber of $p|_{X}$ that is an abelian subvariety, then $p|_X$ is an
abelian subscheme.

\begin{Lemma} \label{lem:abeliansubscheme2}
Let $U$ be an integral noetherian scheme
and let $p : \mathcal{A} \to U$ be an abelian scheme
with zero-section $0_{p}$. 
Let $X$ be a subscheme of $\mathcal{A}$
such that the restriction $p |_{X} : 
X \to U$
is proper and 
smooth
and such that $0_{p}$ factors through ${X} \subset \mathcal{A}$.
Suppose that
there exists a point $s \in U$
such that
the geometric fiber ${X}_{\overline{s}}$ is an
abelian subvariety of $\mathcal{A}_{\overline{s}}$.
Then $p|_{X}$ is an abelian subscheme of $p$.
\end{Lemma}

\Proof
Since 
$p|_X$ is proper and flat and
$X_{\overline{s}}$ is connected,
it follows from
\cite[Proposition~15.5.9~(ii)]{ega4-3}
that
any fiber of $p|_X$ is geometrically connected.
Since $p|_X$ is smooth,
this means that
any fiber of $p|_X$ is geometrically integral.

Let $\alpha : \mathcal{A} \times_U \mathcal{A} \to \mathcal{A}$
denote the addition of $p$.
Set $Z := \alpha ( X \times_U X)$,
the scheme-theoretic image.
Note that $Z$ is an integral scheme.
Indeed,
by \cite[Proposition~4.5.7]{ega4-2},
it follows that
$X \times_U X$ is irreducible.
Further,
since $X \times_U X \to U$ is smooth and $U$ is reduced,
$X \times_U X$ is reduced by
\cite[Proposition~17.5.7]{ega4-4}.
It follows that $X \times_U X$ is integral,
and hence $Z$ is integral.

It follows from
$0_{p} (U) \subset X$
that
$Z \supset X$,
and hence $Z_{u} \supset X_u$ for any $u \in U$. 
Further, 
we have
$Z_{u} = \alpha ( X \times_U X)_{u} =
X_{u} + X_{u}$ as sets.
For any $u \in U$,
since $X_u$ is geometrically irreducible,
it follows that
$Z_u$ is irreducible.

By Lemma~\ref{lem:abeliansubscheme},
we only have to prove $Z \subset X$
as subschemes.
Since $X_{\overline{s}}$ is an abelian variety,
we have
$Z_{\overline{s}} = X_{\overline{s}}$ as sets.
By
the upper semicontinuity of the dimension of fibers
(cf. \cite[Corollaire~13.1.5]{ega4-3}),
there exists
an open neighborhood $V \subset U$ of $s$ such that
$\dim Z_{u} = \dim X_{u}$ for any $u \in V$.
Since
$Z_u$ is irreducible and $Z_u \supset X_u$,
it follows that
$Z_u = X_u$ as sets for $u \in V$.
Therefore,
$Z \times_U V = X \times_U V$ as sets.
Since $Z$ is integral and $Z \supset X$,
it follows that
$Z \times_U V = X \times_U V$ as schemes in fact.
Further, since $X$ is a closed subscheme of $\mathcal{A}$ 
and $Z$ is an integral scheme,
taking the closures of the both sides of this equality
concludes $Z \subset X$ as subschemes.
Thus we 
obtain the lemma.
\QED

The following lemma will be used in the proof of Proposition~\ref{prop:key1}.

\begin{Lemma} \label{lem:generictranslation}
Let $Y$ be an integral noetherian scheme,
let $p : \mathcal{A} \to Y$ be an abelian scheme,
and
let $X$ be a subscheme of $\mathcal{A}$.
Assume
that the restriction $p |_{X} : 
X \to Y$
is proper and surjective
and that the geometric generic fiber $X_{\overline{\eta}}$
is reduced,
where ${\eta}$ is the generic point of $Y$.
Suppose that
there exists a dense subset $S$ of $Y$ such that,
for any $s \in S$,
the subscheme $\left( {X}_{\overline{s}} \right)_{\red}
\subset \mathcal{A}_{\overline{s}}$,
the geometric fiber with its induced reduced subscheme structure, is a
translate of an 
abelian subvariety.
Then,
$X_{\overline{\eta}}$ 
is a translate of an abelian subvariety of $\mathcal{A}_{\overline{\eta}}$.
\end{Lemma}

\Proof
We first claim that there exists
a generically finite dominant morphism $\phi : U \to Y$ 
with $U$ integral
such that, 
if 
$p' : 
U \times_Y \mathcal{A} \to U$
is the base-change of $p$ by $\phi$
and if $Z := U \times_Y X$,
then
the restriction
$p'|_{Z} : Z \to U$
is a proper smooth morphism with
a section $\sigma$.
Indeed,
let $X_\eta$ be the generic fiber of $p|_X$
and let $x \in X_\eta$ be a closed point.
Further,
let $V \subset X$ be the closure of $x$ in $X$.
Then the restriction $V \to Y$ of $p|_{X}$
is a generically finite surjective morphism
such that the base-change $V \times_Y X \to V$ 
of $p|_X$
has a section.
Since $V \times_Y X \to V$
is generically flat,
restricting $V \to Y$ to a dense open subset $U$ of $V$
gives us a generically finite morphism
$\phi : U \to Y$
such that
$p'|_{Z} : Z \to U$
is flat with a section $\sigma$.
Since $p'|_{Z}$ has geometrically reduced generic fiber
by assumption,
\cite[Th\'eor\`eme~12.2.4]{ega4-3} allows us 
moreover to take $\phi$
so that 
any geometric fiber of $p'|_{Z}$ is reduced.
It follows that,
for any $u \in \phi^{-1} (S)$,
we have
$(p'|_Z)^{-1} (u) = ((p'|_Z)^{-1} (u))_{\red}$,
which is smooth by assumption.
This implies that
$p'|_Z$ is smooth over
any $u \in \phi^{-1} (S)$
(cf. \cite[Th\'eor\`eme~17.5.1]{ega4-3}).
Since
$\phi^{-1} (S) \neq \emptyset$
by assumption,
$p'|_Z$ is smooth over an open dense subset of $U$.
Replacing $U$ with this open subset,
we obtain
$p'|_Z :Z \to U$ that is smooth (cf. \cite[Th\'eor\`eme~12.2.4]{ega4-3}).
Thus,
$\phi$ is taken so that
$p'|_Z : Z \to U$ 
is proper and smooth.

Put $Z' := Z - \sigma ( U )$,
the translate of $Z$ by the section $- \sigma ( U )$.
Then $0_{p'} (U) \subset Z'$,
where $0_{p'}$ is the 
zero-section of the abelian scheme $p'$.
Since $p'|_Z : Z \to U$ is smooth and $U \supset \phi^{-1} (S) \neq \emptyset$,
it follows from the assumption 
that
there exits a geometric point $\overline{s}$ 
of $U$ such that
$(p'|_{Z'})^{-1} (\overline{s})$ is an abelian subvariety.
Now, applying Lemma~\ref{lem:abeliansubscheme2}
to these $p' : U \times_Y \mathcal{A}' \to U$ and $Z'$,
we find that $p'|_{Z'}$ is an abelian subscheme of $p'$.
In particular,
the geometric generic fiber of $p'|_{Z'}$
is an abelian subvariety,
and hence the geometric generic fiber $Z_{\overline{\xi}}$
of $p'|_{Z}$
is a translate of an abelian subvariety,
where $\overline{\xi}$ denotes the geometric generic point of $U$.
Since 
$X_{\overline{\eta}} = Z_{\overline{\xi}}$, 
this concludes that $X_{\overline{\eta}}$ is a
translate of an abelian subvariety of $\mathcal{A}_{\overline{\eta}}$.
\QED

\subsection{Relative height} \label{subsec:relativeheight}

We begin by fixing the notation for this subsection.
Let $A$ be an abelian variety over $\overline{K}$
with an even ample line bundle $L$.
Let $\widetilde{Y}$ be an irreducible variety over $k$
and put $Y := \widetilde{Y} \otimes_k \overline{K}$.
Let $p : Y \times {A} \to Y$ be the first projection,
which is an abelian scheme.
Let $X$ be an integral closed subscheme of $Y \times A$
with $p (X ) = Y$.

In this subsection,
we are concerned with a sufficient condition
for $X$ to be of form
$Y \times T$ for some torsion subvariety
$T$ of $A$.
We will give such a condition
in terms of a ``relative height'' function,
which will be denoted by $\mathbf{h}_{X/Y}^{L}$.

To define the relative height function,
we assign
to each point $\tilde{y} \in \widetilde{Y}$
a geometric point $\overline{\tilde{y}_{K}}$
of 
$Y$ 
in the following way.
Taking an algebraic closure
$\mathfrak{k} 
= \overline{k ( \tilde{y})}
$
of the residue field $k ( \tilde{y})$ of $\tilde{y}$
gives rise to a geometric point $\overline{\tilde{y}} :
\Spec \left( 
\mathfrak{k} 
\right) \to \widetilde{Y}$.
The
base-change 
$\left\{ \overline{\tilde{y}} \right\} \times_{\widetilde{Y}}
\left( \widetilde{Y} \times \mathfrak{B} \right)
\cong
\mathfrak{k}
\otimes_{k} \mathfrak{B}$,
denoted by
$\left\{ \overline{\tilde{y}} \right\} \times
\mathfrak{B}$ simply,
is a normal projective variety
over $
\mathfrak{k}
$.
Let 
$\tilde{y}_{K}$
be 
the 
generic point of $\left\{ \overline{\tilde{y}} \right\} \times
\mathfrak{B}$.
We remark that the residue field
$
\mathfrak{k} (\tilde{y}_K)
$
is the
function field of
$\left\{ \overline{\tilde{y}} \right\} \times
\mathfrak{B}$.
Let
$\overline{\mathfrak{k} (\tilde{y}_K)}$
be 
an algebraic closure of $\mathfrak{k} (\tilde{y}_K)$
containing $\overline{K}$.
Then,
the field extensions $\overline{\mathfrak{k} (\tilde{y}_K)} / \mathfrak{k}$
and $\overline{\mathfrak{k} (\tilde{y}_K)} / \overline{K}$
give rise to a morphism
$\Spec \left( \overline{\mathfrak{k} (\tilde{y}_K)} \right)
\to \left\{ \overline{\tilde{y}} \right\} \otimes_{k}
\overline{K}$
and hence
to a geometric point 
\[
\overline{\tilde{y}_K}
:
\Spec
\left( \overline{\mathfrak{k} (\tilde{y}_K)} \right) 
\to 
\widetilde{Y} \otimes_{k} \overline{K}
=Y.
\]
That is the definition of 
$\overline{\tilde{y}_K}$.

\begin{Remark} \label{rem:genericpoint}
If $\tilde{y} = \tilde{\eta}$
is the generic point of $\widetilde{Y}$,
then the corresponding point $\overline{\tilde{\eta}_K}$
is the geometric generic point of $Y$.
\end{Remark}

Let $\OO_{\widetilde{Y}} \otimes_k \mathcal{H}$ be the pull-back
of $\mathcal{H}$ by the natural projection 
$\widetilde{Y} \times_{\Spec (k)} \mathfrak{B} \to \mathfrak{B}$.
Since
$\mathfrak{k} (\tilde{y}_K)$ is the function field of 
the normal projective variety 
$\left\{ \overline{\tilde{y}} \right\} \times
\mathfrak{B}$
equipped with an ample line bundle $\mathfrak{k} \otimes_{k} \mathcal{H}
= \left( \OO_{\widetilde{Y}} \otimes_k \mathcal{H} \right)|_{\left\{ \overline{\tilde{y}} 
\right\} \times \mathfrak{B}}$,
the notion of height
over the base field $\mathfrak{k} (\tilde{y}_K)$ makes sense.
Indeed,
we
consider the
abelian variety
$
\left\{ \overline{\tilde{y}_K} \right\} \times_{Y}
(Y \times A)
\cong
\overline{\mathfrak{k} (\tilde{y}_K)}
\otimes_{\overline{K}} 
A
$
over $\overline{\mathfrak{k} (\tilde{y}_K)}$,
which is denoted by $
\left\{ \overline{\tilde{y}_K} \right\} \times
A$.
This abelian variety is equipped with the even ample line bundle
$
\overline{\mathfrak{k} (\tilde{y}_K)}
\otimes_{\overline{K}} 
L
= \left( \OO_Y \otimes_{\overline{K}} L \right)|_{
\left\{ \overline{\tilde{y}_K} \right\} \times
A}$,
where $\OO_Y \otimes_{\overline{K}} L$ is the pull-back of $L$
by the natural projection $Y \times A \to A$.
Since
the fiber
$X_{\overline{\tilde{y}_K}}$
of
$p|_{X} : X \to Y$ over $\overline{\tilde{y}_K}$
is a closed subscheme of 
$
\left\{ \overline{\tilde{y}_K} \right\} \times
A$,
we can therefore consider the canonical height of $X_{\overline{\tilde{y}_K}}$
with respect to this line bundle
if $X_{\overline{\tilde{y}_K}}$
has pure dimension.

Now,
varying $\tilde{y}$
defines
a function $\mathbf{h}^{L}_{X/Y}
$
which assigns to each $\tilde{y}$ the height of 
$X_{\overline{\tilde{y}_K}}$ with respect to 
$
\overline{\mathfrak{k} (\tilde{y}_K)}
\otimes_{\overline{K}} 
L
= \left( \OO_Y \otimes_{\overline{K}} L \right)|_{
\left\{ \overline{\tilde{y}_K} \right\} \times
A}$.
To be precise,
let $\widetilde{Y}_{\mathrm{pd}}$
be a subset
of $\widetilde{Y}$ given by
\[
\widetilde{Y}_{\mathrm{pd}}
:=
\left\{
\left.
\tilde{y} \in \widetilde{Y}
\right|
\text{
$X_{\overline{\tilde{y}_K}}$
has pure dimension
$\dim (X) - \dim (Y)$}
\right\}
.
\]
Then we
define a function
$\mathbf{h}^{L}_{X/Y}$ on $\widetilde{Y}_{\mathrm{pd}}$
by
\addtocounter{Claim}{1}
\begin{align} \label{eq:defofmathbfh}
\mathbf{h}^{L}_{X/Y} ( \tilde{y} )
:=
\widehat{h}_{
\overline{\mathfrak{k} (\tilde{y}_K)}
\otimes_{\overline{K}} 
L
} 
\left(
X_{\overline{\tilde{y}_K}}
\right)
=
\widehat{h}_{\left( \OO_Y \otimes_{\overline{K}} L \right)|_{
\left\{ \overline{\tilde{y}_K} \right\} \times
A}} 
\left(
X_{\overline{\tilde{y}_K}}
\right)
.
\end{align}

\begin{Remark} \label{rem:mathbfh-nonnegative}
For any $\tilde{y} \in \widetilde{Y}_{\mathrm{pd}}$
and for any irreducible component $Z$ of $X_{\overline{\tilde{y}_K}}$,
we have 
\[
0 \leq 
\widehat{h}_{
\overline{\mathfrak{k} (\tilde{y}_K)}
\otimes_{\overline{K}} 
L
} 
(
Z
)
\leq
\widehat{h}_{
\overline{\mathfrak{k} (\tilde{y}_K)}
\otimes_{\overline{K}} 
L} 
\left(
X_{\overline{\tilde{y}_K}}
\right)
=
\mathbf{h}^{L}_{X/Y} (\tilde{y}),
\]
which follows from 
the non-negativity
of the canonical height
\cite[Theorem~11.18~(e)]{gubler0}
and the linearity of the canonical height on cycles
\cite[Theorem~3.5~(a)]{gubler3}.
\end{Remark}

The generic point $\tilde{\eta}$ of 
$\widetilde{Y}$ sits in $\widetilde{Y}_{\mathrm{pd}}$.
Indeed, 
since $X$ is irreducible,
the generic fiber
of $X \to Y$ is also irreducible
and has dimension $\dim (X) - \dim (Y)$.
Therefore the geometric generic fiber of $X \to Y$
has pure dimension $\dim (X) - \dim (Y)$.
Since
the point $\overline{\tilde{\eta}_{K}}$
associated to $\tilde{\eta}$
is the geometric generic point of $Y$ (cf. Remark~\ref{rem:genericpoint}),
we have $\tilde{\eta} \in \widetilde{Y}_{\mathrm{pd}}$.

The following proposition gives us a condition
in terms of the relative height
for $X$ to be 
the product of $Y$ with a torsion subvariety of $A$.

\begin{Proposition} \label{prop:key1}
Let $A$ be an abelian variety over $\overline{K}$
with an even ample line bundle $L$.
Let $\widetilde{Y}$ be an irreducible variety over $k$
and put $Y := \widetilde{Y} \otimes_k \overline{K}$.
Let $p : Y \times {A} \to Y$ be the first projection.
Let $X$ be an integral closed subscheme of $Y \times {A}$
with $p (X ) = Y$.
Assume that $A$ 
has trivial $\overline{K}/k$-trace.
Suppose that the following conditions are satisfied.
\begin{enumerate}
\renewcommand{\labelenumi}{(\alph{enumi})}
\item
We have
$\mathbf{h}^{L}_{X/Y} (\tilde{\eta}) = 0$,
where $\tilde{\eta}$ is the generic point
of $\widetilde{Y}$.
\item
There exists a dense subset $S \subset Y$
such that,
for any  
$s \in S$,
each irreducible component of
the geometric fiber
$p|_X^{-1} (\overline{s})$
with its induced reduced subscheme structure
is a
torsion subvariety of $\{ \overline{s} \} \times A$.
\end{enumerate}
Then, there exists a torsion subvariety $T \subset A$
such that $X = Y \times T$.
\end{Proposition}

\Proof
We use the following notation and convention
in this proof.
For a morphism $\phi : U \to Y$,
put $X_U := U \times_Y X$
and let $(p|_{X})_U : X_U \to U$ denote the base-change of $p|_{X}$
by $\phi$.
When we say an irreducible component $Z$ of a scheme,
we regard $Z$ as a closed subscheme
\emph{with
its induced reduced subscheme structure}.

There exists a generically finite dominant morphism
$\phi : U \to Y$ with $U$ integral
and an
irreducible component $Z$ of $X_U$
such that
the restriction $q : Z \to U$ 
of the morphism $(p|_{X})_{U} : X_U \to U$ to $Z$
is a proper flat morphism with geometrically integral fibers.
Indeed,
there exists a generically finite dominant morphism $\phi : U \to Y$ 
with $U$ integral
such that,
for any irreducible component $Z$ of $X_U$,
the restriction 
$(p|_{X})_{U}$ to $Z$
is a surjective morphism with
geometrically integral generic fiber.
Shrinking $U$ if necessary,
we may take $\phi$ 
so that
$Z \to U$
is flat,
and
by \cite[Th\'eor\`eme~12.2.4]{ega4-3},
we may take it moreover so that any fiber of $Z \to  U$
is geometrically integral.
The properness of $(p|_{X})_{U}$
follows from the fact that X is proper over Y.

Let 
$\mathfrak{k}$ be an
algebraic closure of
the residue field at $\tilde{\eta}$
and 
let
$\overline{\tilde{\eta}_K} : \Spec \left( 
\overline{\mathfrak{k} (\tilde{\eta}_K)} \right) \to Y$ 
be the corresponding geometric point of $Y$
defined at the beginning of this subsection.
Let $\overline{\xi}$ be the geometric generic point of $U$.
Then $\phi \left( \overline{\xi} \right)$ is the 
geometric generic point of $Y$,
and
it follows from
Remark~\ref{rem:genericpoint} 
that $\phi \left( \overline{\xi} \right) = \overline{\tilde{\eta}_K}$.

We apply Lemma~\ref{lem:generictranslation}
to $U$, $U \times_{\Spec \left( \overline{K} \right)} A \to U$ and $Z$
in place of $Y$, $\mathcal{A} \to Y$ and $X$.
To do that,
we remark that, 
$\phi^{-1} (S)$ is dense in $U$ by the assumption (b).
Further, we remark
that
for any $u \in U$,
$Z_{\overline{u}} := q^{-1} (\overline{u})$
is an integral subscheme
that appears as
an irreducible component
of $(p|_{X})_{U}^{-1} (\overline{u}) = (p|_{X})^{-1} ( \phi ( \overline{u} ))$.
By assumption,
it follows
that, for
any $s \in \phi^{-1}(S)$,
$Z_{\overline{s}}$ is a torsion subvariety.
Then,
Lemma~\ref{lem:generictranslation}
concludes that
there exist a point $\sigma \in A \left(
\overline{\mathfrak{k} (\tilde{\eta}_K)}
\right)$
and an abelian subvariety $G' \subset \left\{ \overline{\xi} \right\} 
\times
A$ such that
$Z_{\overline{\xi}} = G' + \sigma$.

We prove that 
$Z_{\overline{\xi}}$ is a torsion subvariety.
By
the assumption (a),
we have
\[
\widehat{h}_{
\overline{\mathfrak{k} (\tilde{y}_K)}
\otimes_{\overline{K}} 
L
}
\left( \left( p|_{X} \right)^{-1}
\left( \phi \left( \overline{\xi} \right) \right)  \right)
=
\widehat{h}_{
\overline{\mathfrak{k} (\tilde{y}_K)}
\otimes_{\overline{K}} 
L
}
\left(  X_{\overline{\tilde{\eta}_K}}  \right)
=
\mathbf{h}^{L}_{X/Y} ( \tilde{\eta} )
=
0
.
\]
Since $Z_{\overline{\xi}}$
is an irreducible component of 
$\left( p|_{X} \right)^{-1}( \phi(\overline{\xi})) = X_{\overline{\tilde{\eta}_K}}$,
we then see from
Remark~\ref{rem:mathbfh-nonnegative}
that
$\widehat{h}_{
\overline{\mathfrak{k} (\tilde{y}_K)}
\otimes_{\overline{K}} 
L
}
\left( Z_{\overline{\xi}} \right) = 0$.
By
Proposition~\ref{prop:dense-height0},
it follows that
$Z_{\overline{\xi}}$ has dense small points,
and hence
$Z_{\overline{\xi}} / G' \subset
\left(
\left\{ \overline{\xi} \right\} \times A
\right)  / G'$ has dense small points
by \cite[Lemma~2.1]{yamaki5}.
This means that the image of $\sigma$ in
$\left(
\left\{ \overline{\xi} \right\} \times A
\right)  / G'$
has canonical height $0$.
Since 
$A$ has trivial 
$\overline{K}/k$-trace, 
$\left\{ \overline{\xi} \right\} \times A$
has trivial $
\overline{\mathfrak{k} (\tilde{\eta}_K)}/
\mathfrak{k}$-trace
by Lemma~\ref{lem:trivial-trivial},
which will be proved in the appendix.
Hence,
by \cite[Lemma~1.5]{yamaki5},
$\left(
\left\{ \overline{\xi} \right\} \times A
\right)  / G'$
has trivial $
\overline{\mathfrak{k} (\tilde{\eta}_K)}/
\mathfrak{k}$-trace.
By \cite[Theorem~5.4]{lang2},
it follows that
the image of 
$\sigma$ in
$\left(
\left\{ \overline{\xi} \right\} \times A
\right)  / G'$
is a torsion point.
Since a surjective homomorphism between abelian varieties
over an algebraically closed field
induces a surjective homomorphism between the
groups of torsion points,
there exists a torsion point $\tau \in
\left\{ \overline{\xi} \right\} \times A
$
having the same image 
as $\sigma$
in 
$\left(
\left\{ \overline{\xi} \right\} \times A
\right)  / G'$.
Then $Z_{\overline{\xi}} = G' + \sigma = G ' + \tau$,
which shows that $Z_{\overline{\xi}}$
is a torsion subvariety of 
$\left\{ \overline{\xi} \right\} \times A
$.

By
Proposition~\ref{prop:Chow},
therefore,
there exists
a torsion subvariety $T$ of $A$
such that $Z_{\overline{\xi}} = \left\{ \overline{\xi} \right\} \times T$.
We look at the closed subschemes $Y \times T$ and $X$ of $Y \times A$.
We remark that,
since $\phi \left( \overline{\xi} \right) = \overline{\tilde{\eta}_{K}}$,
$Z_{\overline{\xi}}$
equals the geometric generic fiber of
$p|_{Y \times T} :  Y \times T \to Y$.
Since $Z_{\overline{\xi}}$ is 
a closed subscheme
of 
$X_{\overline{\tilde{\eta}_{K}}}$,
it follows that
the generic fiber of $p|_{Y \times T}$
is
a closed subscheme of the generic fiber of $p|_X : X \to Y$.
Since $Y \times T$ is integral,
this implies that $Y \times T \subset X$ as subschemes.
Since $X$ is integral and $\dim (Y \times T) = \dim Z = \dim X$,
we then conclude that $Y \times T = X$ as subschemes.
Thus we obtain the proposition.
\QED

\subsection{Generic constantness of relative heights}

The purpose of this subsection is to show that
the relative height
$\mathbf{h}^{L}_{X/Y}$ is generically constant
when $A$ is nowhere degenerate.

We begin with a technical lemma.

\begin{Lemma} \label{lem:flatness1}
Let $\widetilde{Y}$ be a proper irreducible variety over $k$.
Let $K'$ be a finite extension 
of $K$
and let $\mathfrak{B}'$
be the normalization of $\mathfrak{B}$ in $K'$.
Let
$h : \mathscr{X} \to \widetilde{Y} \times_{\Spec (k)} \mathfrak{B}'$
be a proper surjective morphism with $\mathscr{X}$ integral.
Then, there exists
a dense open subset $\widetilde{V} \subset \widetilde{Y}$
satisfying the following conditions.
\begin{enumerate}
\renewcommand{\labelenumi}{(\alph{enumi})}
\item
Let $\pr_{\widetilde{Y}} : \widetilde{Y} \times_{\Spec (k)} \mathfrak{B}'
\to \widetilde{Y}$ be the first projection.
Then 
the restriction
$\widetilde{V} \times_{\widetilde{Y}} \mathscr{X}
\to \widetilde{V}$
of $\pr_{\widetilde{Y} } \circ h
:
\mathscr{X} \to \widetilde{Y}$
is flat.
\item
For any $\tilde{y} \in \widetilde{V}$,
the restriction
$h|_{\{ \tilde{y} \} \times_{\widetilde{Y}} \mathscr{X}}
: 
\{ \tilde{y} \} \times_{\widetilde{Y}} \mathscr{X}
\to \{ \tilde{y} \} \times_{\Spec (k)} \mathfrak{B}'
$
of $h$ is flat over 
any point of codimension $1$
in $\{ \tilde{y} \} \times_{\Spec (k)} \mathfrak{B}'$.
\end{enumerate}
\end{Lemma}

\Proof
First,
since $\mathscr{X}$ is an integral scheme,
the generic flatness gives us
a dense open subset $\widetilde{V}_1$ of $\widetilde{Y}$
such that the
restriction
$\widetilde{V}_1 \times_{\widetilde{Y}} \mathscr{X}
\to \widetilde{V}_1$
of $\pr_{\widetilde{Y}} \circ h : \mathscr{X} \to \widetilde{Y}$
is flat.

Next,
we construct $\widetilde{V}$ which will suffice (b).
Let $\widetilde{Y}_{\mathrm{reg}}$
be the subset of regular points
of $\widetilde{Y}$.
It is a dense open subset.
Since $\mathscr{X}$ is an 
integral scheme and $\widetilde{Y}_{\mathrm{reg}} \times_{\Spec (k)} \mathfrak{B}'$
is normal,
there exists a closed subset $\mathcal{Z}$ of 
$\widetilde{Y}_{\mathrm{reg}} \times_{\Spec (k)} \mathfrak{B}'$ with
$\codim \left( 
\mathcal{Z} , \widetilde{Y}_{\mathrm{reg}} \times_{\Spec (k)} \mathfrak{B}'
\right) 
\geq 2$ such that
$h : \mathscr{X} \to \widetilde{Y} \times_{\Spec (k)} \mathfrak{B}'$ 
is flat over $
\left(
\widetilde{Y}_{\mathrm{reg}} \times_{\Spec (k)} \mathfrak{B}' 
\right) \setminus \mathcal{Z}$.
Then,
for any $\tilde{y} \in \widetilde{Y}_{\mathrm{reg}}$,
$\{ \tilde{y} \}
\times_{\widetilde{Y}}
\mathcal{Z}$ is a closed subset
of $\{ \tilde{y} \}
\times_{\Spec (k)} \mathfrak{B}'$,
and
the base-change
$\{ \tilde{y} \} \times_{\widetilde{Y}} \mathscr{X}
\to \{ \tilde{y} \}
\times_{\Spec (k)} \mathfrak{B}'$
of $h$
is flat over 
$\left( \{ \tilde{y} \}
\times_{\Spec (k)} \mathfrak{B}' \right)
\setminus
\left( \{ \tilde{y} \} \times_{\widetilde{Y} } \mathcal{Z} \right)$.
Let $\tilde{\eta}$ be the generic point of $\widetilde{Y}$.
Then
$\tilde{\eta} \in \widetilde{Y}_{\mathrm{reg}}$,
and
$\codim 
\left( 
\{ \tilde{\eta} \}
\times_{\widetilde{Y}}
\mathcal{Z} ,
\{ \tilde{\eta} \}
\times_{\Spec (k)} \mathfrak{B}'
\right) \geq 2$.
It follows that there exists a dense open subset $\widetilde{V}_2
\subset \widetilde{Y}_{\mathrm{reg}}$ such that,
for any $\tilde{y} \in \widetilde{V}_2$,
we have
$\codim 
\left( 
\{ \tilde{y} \} \times_{\widetilde{Y} } \mathcal{Z},
\{ \tilde{y} \}
\times_{\Spec (k)} \mathfrak{B}'
\right) \geq 2$.
This means that 
the condition (b) is satisfied for
$V = V_2$.

Finally, we set $V := V_1 \cap V_2$.
Then it is a dense open subset of $\widetilde{Y}$ 
which satisfies both the conditions.
Thus we conclude the lemma.
\QED

\begin{Remark} \label{rem:flatness2}
Under the setting of Lemma~\ref{lem:flatness1},
let $\widetilde{V} \subset \widetilde{Y}$ be an open subset 
as
in the lemma.
Then the above proof shows that 
there exists a closed subset $\mathcal{Z}'
\subset \widetilde{V} \times_{\Spec (k)} \mathfrak{B}'$
such that
$\codim 
\left( 
\{ \tilde{y} \} \times_{\widetilde{Y} } \mathcal{Z}',
\{ \tilde{y} \}
\times_{\Spec (k)} \mathfrak{B}'
\right) \geq 2$
for any $\tilde{y} \in \widetilde{V}$
and that
$h
: \mathscr{X} \to \widetilde{Y} \times_{\Spec (k)} \mathfrak{B}'$ is flat over $
\left( \widetilde{V} \times_{\Spec (k)} \mathfrak{B}'
\right)
\setminus 
\mathcal{Z}'$.
Indeed, if $\mathcal{Z}$ is the closed subset of $\widetilde{Y}_{\mathrm{reg}} \times_{\Spec (k)} \mathfrak{B}'$ as in the proof,
then  
$\mathcal{Z}' :=
\mathcal{Z} \cap 
\left( 
\widetilde{V} \times_{\Spec (k)} \mathfrak{B}' \right)$
has those properties.
\end{Remark}

Now, we show the generic constantness of the relative height:

\begin{Proposition} \label{prop:hgenericallyconstant}
Let $A$ be an abelian variety over $\overline{K}$
with an even ample line bundle $L$.
Let $\widetilde{Y}$ be an irreducible variety over $k$
and put $Y := \widetilde{Y} \otimes_k \overline{K}$.
Let $p : Y \times A \to Y$ be the first projection
and
let $X$ be an integral closed subscheme of $Y \times {A}$
with $p (X ) = Y$.
Assume that $A$ is nowhere degenerate.
Then
there exists a dense open subset $\widetilde{V}$
of $\widetilde{Y}$ contained in $\widetilde{Y}_{\mathrm{pd}}$ such that
$\mathbf{h}^{L}_{X/Y}$
is constant over $\widetilde{V}$.
\end{Proposition}

\Proof
Let
$K'$, $\mathfrak{B}'$, $f: \mathscr{A} \to \mathfrak{B}'$
and $\mathscr{L}$ be those as in 
Proposition~\ref{prop:model1} for $A$ and $L$.
Taking a finite extension of $K'$ if necessary,
we
assume that $X$ can be defined over $K'$.
Let
$\mathscr{X}$ be the closure of $X$
in $\widetilde{Y} \times_{\Spec (k)} \mathscr{A}$.
Since $X$ is integral, $\mathscr{X}$ is also integral.
Let 
$h : \mathscr{X} \to \widetilde{Y} \times_{\Spec (k)} \mathfrak{B}'$
be the restriction to $\mathscr{X}$
of
the base-change
$f_{\widetilde{Y}} : \widetilde{Y} \times_{\Spec (k)} \mathscr{A} \to
\widetilde{Y} \times_{\Spec (k)} \mathfrak{B}'$
of $f$.
Since
$f_{\widetilde{Y}}$
is proper,
$h$ is also proper.
The restriction of $h$ to the geometric generic fiber over $\mathfrak{B}'$
equals
$p|_{X}$.
Since 
$p|_{X}$ is surjective and $\mathfrak{B}'$
is irreducible, it follows 
that $h$ is surjective.
By Lemma~\ref{lem:flatness1},
we then take a dense open subset $\widetilde{V}$ satisfying the conditions
(a) and (b)
of Lemma~\ref{lem:flatness1}.

Let us show that
$\widetilde{V} \subset \widetilde{Y}_{\mathrm{pd}}$.
Let $\tilde{y} \in \widetilde{V}$ be any point.
Recall from \S~\ref{subsec:relativeheight}
that
$\overline{\tilde{y}_K}$ is the geometric generic point
of $\{ {\tilde{y}} \} \times_{\Spec (k)} \mathfrak{B}'$
and hence a geometric point of $\widetilde{Y} \times_{\Spec (k)} \mathfrak{B}'$,
and furthermore recall that
$X_{\overline{\tilde{y}_K}}$ is nothing but the fiber
of $h : \mathscr{X} \to \widetilde{Y} \times_{\Spec (k)} \mathfrak{B}'$
over $\overline{\tilde{y}_K}$.
Let
$\mathcal{Z}'$ be a closed subset of 
$\widetilde{V} \times_{\Spec (k)} \mathfrak{B}'$
as in Remark~\ref{rem:flatness2}.
Since
 $\left(
\left\{ \tilde{y} 
\right\} \times_{\Spec (k)} \mathfrak{B}'
\right)
\setminus
\left( \{ \tilde{y} \} \times_{\widetilde{Y}}
\mathcal{Z}' \right)$
is a dense open subset of 
$\left\{ \tilde{y} 
\right\} \times_{\Spec (k)} \mathfrak{B}'$
and $\overline{\tilde{y}_K}$
is the geometric \emph{generic} point of $\{ \tilde{y} \} 
\times_{\Spec (k)} \mathfrak{B}'$,
$\overline{\tilde{y}_K}$ is also
the geometric generic point of 
$\left( \{ \tilde{y} \} 
\times_{\Spec (k)} \mathfrak{B}' \right)
\setminus
\left( \{ \tilde{y} \} \times_{\widetilde{Y}}
\mathcal{Z}' \right)$
and hence
is a geometric point
of
$\left(
\widetilde{V} \times_{\Spec (k)} \mathfrak{B}'
\right)
\setminus
\mathcal{Z}'$.
%
Since $h$ is flat over $\left(
\widetilde{V} \times_{\Spec (k)} \mathfrak{B}'
\right)
\setminus
\mathcal{Z}'$
and $\mathscr{X}$ is irreducible,
it follows that
the fiber $X_{\overline{\tilde{y}_K}}$ 
of $h$ over ${\overline{\tilde{y}_K}}$
has pure dimension 
\[
\dim (\mathscr{X}) - \dim 
\left(  \widetilde{V} \times_{\Spec (k)} \mathfrak{B}' \right)
%
=
\dim (X) - \dim (Y).
\]
This shows $\widetilde{y} \in \widetilde{Y}_{\mathrm{pd}}$,
and thus $\widetilde{V} \subset \widetilde{Y}_{\mathrm{pd}}$.

We are going to prove that
$\mathbf{h}^{L}_{X/Y}
$
is constant over $\widetilde{V}$.
To do that,
we 
describe
 $\mathbf{h}^{L}_{X/Y}
$ in terms of intersection
product on models.
For each $\tilde{y} \in \widetilde{V}$,
let $\overline{\tilde{y}} : \Spec \left( \overline{k(\tilde{y})} \right)
\to \widetilde{Y}$ be the geometric point arising from $\tilde{y}$.
Since
the point $\overline{\tilde{y}_K}$
is
the geometric generic point of $\left\{ \overline{\tilde{y}} \right\} \times
\mathfrak{B}'$,
we note that
$\left\{ \overline{\tilde{y}_K} \right\} \times A$
is
the geometric generic fiber of
$
f_{\overline{\tilde{y}}}
:
\left\{ \overline{\tilde{y}} \right\} \times
\mathscr{A}
\to
\left\{ \overline{\tilde{y}} \right\} \times
\mathfrak{B}'
$,
where $f_{\overline{\tilde{y}}}$ is the 
restriction of $f_{\widetilde{Y}}$.
Let 
$\OO_{\widetilde{Y}} \otimes_k \mathscr{L}$ be the pull-back
of $\mathscr{L}$ by the canonical projection 
$\widetilde{Y } \times_{\Spec (k)} \mathscr{A} \to \mathscr{A}$
and let $\OO_Y \otimes_{\overline{K}} L$ be the pull-back
of $L$ by the canonical projection $Y \times A \to A$.
Then we see that $\left( \left\{ \overline{\tilde{y}} \right\} \times
\mathscr{A} , \left(
\OO_{\widetilde{Y}} \otimes_k
\mathscr{L}
\right)
|_{\left\{ \overline{\tilde{y}} \right\} \times
\mathscr{A}}
\right)$
is a model of $\left(
\left\{ \overline{\tilde{y}_K} \right\} \times A ,
\left( \OO_{Y} \otimes_{\overline{K}} L \right)|_{
\left\{ \overline{y_K} \right\} \times A}
\right)$
satisfying the conditions of Proposition~\ref{prop:model1}.

Let $\pr_{\widetilde{Y}} : \widetilde{Y} \times_{\Spec (k) } \mathfrak{B}'
\to \widetilde{Y}$ be the projection.
For any $\tilde{y} \in \widetilde{V}$,
put $\mathscr{X}_{\overline{\tilde{y}}} :=
\left( \pr_{\widetilde{Y}} \circ h \right)^{-1} (\overline{\tilde{y}})$,
which
is a closed subscheme 
of 
$\left\{ \overline{\tilde{y}} \right\} \times
\mathscr{A}$.
We remark that
$p|_X : X \to Y$ equals the restriction of $h : \mathscr{X} \to 
\widetilde{Y} \times_{\Spec (k)} \mathfrak{B}'$
to the geometric generic fiber over $\mathfrak{B}'$.
Then we see that
$X_{\overline{\tilde{y}_K}} := \left( p|_{X} \right)^{-1} 
\left(
\overline{\tilde{y}_K}
\right)$
is the geometric generic fiber
of $h|_{\mathscr{X}_{\overline{\tilde{y}}}} :
\mathscr{X}_{\overline{\tilde{y}}} \to
\left\{ \overline{\tilde{y}} \right\} \times
\mathfrak{B}'$.
By the condition (b)
of
Lemma~\ref{lem:flatness1}
for $\widetilde{V}$,
the proper morphism $h|_{\mathscr{X}_{\overline{\tilde{y}}}}
$
is 
flat over any 
point of $\left\{ \overline{\tilde{y}} \right\} \times
\mathfrak{B}'$ of codimension $1$.
Further,
since $\mathscr{X}$ is an integral scheme and
since $\pr_{\widetilde{Y}} \circ h$ is flat over $\widetilde{V}$
by the condition (a) of Lemma~\ref{lem:flatness1},
we note that $\mathscr{X}_{\overline{\tilde{y}}}$ is pure dimensional.
Thus
it is a model of $X_{\overline{\tilde{y}_K}}$
in our sense.

Recalling (\ref{eq:defofmathbfh}),
we then apply
Lemma~\ref{lem:height-intersection}
to
obtain
\addtocounter{Claim}{1}
\begin{align} \label{eq:relativeheight-intersection}
\mathbf{h}^{L}_{X/Y}
\left( \tilde{y} \right)
=
\frac{
\deg_{\left( \overline{k(\tilde{y})} \otimes_{k} \mathcal{H}' \right)}
\left( 
f_{\overline{\tilde{y}}}
\right)_{\ast}
\left(
\cherncl_1
\left(
\left(
\OO_{\widetilde{Y}} \otimes_k
\mathscr{L}
\right)
|_{\left\{ \overline{\tilde{y}} \right\} \times
\mathscr{A}}
\right)^{\cdot (\dim (X) - \dim (Y) + 1)}
\cdot
\left[ \mathscr{X}_{\overline{\tilde{y}}} \right]
\right)
}{[K' : K]}
,
\end{align}
where
$\mathcal{H}'$ is the pull-back of $\mathcal{H}$
by the morphism $\mathfrak{B}' \to \mathfrak{B}$
and $\overline{k(\tilde{y})} \otimes_{k} \mathcal{H}' $
is the pull-back of $\mathcal{H}'$
by the natural morphism $\left\{ \overline{\tilde{y}} \right\}
\times \mathfrak{B}' \to \mathfrak{B}'$.
Using the projection formula,
we see 
\begin{align*} 
&\deg_{\left( \overline{k(\tilde{y})} \otimes_{k} \mathcal{H}' \right)}
\left( 
f_{\overline{\tilde{y}}}
\right)_{\ast}
\left(
\cherncl_1
\left(
\left(
\OO_{\widetilde{Y}} \otimes_k
\mathscr{L}
\right
|_{\left\{ \overline{\tilde{y}} \right\} \times
\mathscr{A}}
\right)^{\cdot (\dim (X) - \dim (Y) + 1)}
\cdot
\left[ \mathscr{X}_{\overline{\tilde{y}}} \right]
\right)
\\
&=
\deg
\left(
\cherncl_1
\left(
\left( \OO_{\widetilde{Y}} \otimes_k f^{\ast} ( \mathcal{H}' ) \right)
|_{\left\{ \overline{\tilde{y}} \right\} \times
\mathscr{A}}
\right)^{\cdot (b-1)}
\cdot
\cherncl_1
\left(
\left(
\OO_{\widetilde{Y}} \otimes_k
\mathscr{L}
\right)
|_{\left\{ \overline{\tilde{y}} \right\} \times
\mathscr{A}}
\right)^{\cdot (\dim (X) - \dim (Y) + 1)}
\cdot
\left[ \mathscr{X}_{\overline{\tilde{y}}} \right]
\right)
\\
&=
\deg
\left(
\cherncl_1
\left(
\OO_{\widetilde{Y}} \otimes_k f^{\ast} ( \mathcal{H}' )
\right)^{\cdot (b-1)}
\cdot
\cherncl_1
\left(
\OO_{\widetilde{Y}} \otimes_k
\mathscr{L}
\right)^{\cdot (\dim (X) - \dim (Y) + 1)}
\cdot
\left[ \mathscr{X}_{\overline{{\tilde{y}}}} \right]
\right)
,
\end{align*}
and hence
\addtocounter{Claim}{1}
\begin{align} \label{eq:mathbfh-intersection}
\mathbf{h}^{L}_{X/Y}
\left( \tilde{y} \right)
=
\frac{
\deg
\left(
\cherncl_1
\left(
\OO_{\widetilde{Y}} \otimes_k f^{\ast} ( \mathcal{H}' )
\right)^{\cdot (b-1)}
\cdot
\cherncl_1
\left(
\OO_{\widetilde{Y}} \otimes_k
\mathscr{L}
\right)^{\cdot (\dim (X) - \dim (Y) + 1)}
\cdot
\left[ \mathscr{X}_{\overline{\tilde{y}}} \right]
\right)}{[K' :K]}
.
\end{align}
%
By the condition (a) of Lemma~\ref{lem:flatness1},
$\pr_{\widetilde{Y}} \circ h$
is flat over 
$\widetilde{V}$.
By \cite[Theorem~10.2]{fulton},
it follows that
the intersection number on the left-hand side
in (\ref{eq:mathbfh-intersection})
is independent of $\tilde{y} \in \widetilde{V}$,
and hence
$\mathbf{h}^{L}_{X/Y}$ is also independent of $\tilde{y} \in \widetilde{V}$.
Thus we conclude the proposition.
\QED


\section{Application to the geometric Bogomolov conjecture}
\label{sect:application}

In this section, we prove Theorem~\ref{thm:main2intro}
and
reduce
the geometric Bogomolov conjecture for any abelian varieties
to the conjecture for
those
abelian varieties which are nowhere degenerate and have trivial trace.
We begin with a proposition which
will be the key to this goal.

\begin{Proposition} \label{prop:final}
Let $A$ be a nowhere degenerate abelian variety over $\overline{K}$
and let $L$ be an even ample line bundle on $A$.
Let
$\widetilde{B}$ be an abelian variety over $k$
and set $B := \widetilde{B} \otimes_k \overline{K}$.
Let $X$ be an irreducible closed subvariety of $B \times A$.
Let $\pr_B : B \times A \to B$ be the natural projection
and
set $Y := \pr_B (X)$.
Suppose that $X$ has dense small points.
Then,
there exists a closed subvariety $\widetilde{Y} \subset \widetilde{B}$
such that $Y = \widetilde{Y} \otimes_k \overline{K}$.
Furthermore,
$\mathbf{h}^{L}_{X/Y} (\tilde{y}) = 0$ holds for
general $\tilde{y} \in \widetilde{Y} (k)$
\end{Proposition}

\Proof
Let $\widetilde{M}$
be an even and very ample line bundle on $\widetilde{B}$
and set
$M := \widetilde{M} \otimes_k \overline{K}$.
Since $X$ has dense small points,
so does $Y$
by
\cite[Lemma~7.7]{yamaki6},
and hence $\widehat{h}_{M} (Y) = 0$ by Proposition~\ref{prop:dense-height0}.
It follows from Proposition~\ref{prop:GBCforconstant}
that
there exists
a closed subvariety $\widetilde{Y}$ of $\widetilde{B}$ such that 
$Y = \widetilde{Y} \otimes_k \overline{K}$.
Thus we obtain the first part of the proposition.

Next, we discuss the second part of the proposition.
For $A$ and $L$,
let $K'$,
$f : \mathscr{A} \to \mathfrak{B}'$,
and a line bundle $\mathscr{L}$ on $\mathscr{A}$
be as in Proposition~\ref{prop:model1}.
Taking a finite extension of $K'$ if necessary,
we may assume that $X$ can be defined over $K'$.
Let $\mathcal{H}'$ be the pull-back of $\mathcal{H}$ by the 
morphism $\mathfrak{B}' \to \mathfrak{B}$.

Set $\mathscr{B} := \widetilde{B} \times_{\Spec (k)} \mathfrak{B}'$.
Then the canonical projection $\mathscr{B} \to \mathfrak{B}'$
is a model of $B$,
and
$\widetilde{Y}
\times_{\Spec (k)} \mathfrak{B}'$
equals
the closure of $Y$ in $\mathscr{B}$.
Let $\pr_{\mathscr{A}} : \mathscr{B} \times_{\mathfrak{B}'}
\mathscr{A} \to \mathscr{A}$
be 
the canonical projection
and set
$\varphi := f \circ \pr_{\mathscr{A}}
:
\mathscr{B} \times_{\mathfrak{B}'}
\mathscr{A} \to \mathfrak{B}'
$.
Then
$\varphi$
is a model of
$B \times A$.
Let $\mathscr{X}$ be the closure of $X$ in the model $\mathscr{B}
\times_{\mathfrak{B}'} \mathscr{A}$.
Since $X$ is integral,
$\mathscr{X}$ is also integral.
Let $\pr_{\mathscr{B}} : \mathscr{B} \times_{\mathfrak{B}'} \mathscr{A}
\to \mathscr{B}$
be the natural projection.
Since $\pr_{B}  (X) = Y$,
it follows that
$\pr_{\mathscr{B}} ( \mathscr{X} ) = 
\widetilde{Y}
\times_{\Spec (k)} \mathfrak{B}'$.
%
%
Let
$h :
\mathscr{X} \to \widetilde{Y}
\times_{\Spec (k)} \mathfrak{B}'$
the morphism given from $\pr_{\mathscr{B}}$
by restriction.
Then
$h$ is proper and surjective.
Now,
applying
Lemma~\ref{lem:flatness1} to this $h$,
we obtain a dense open subset
$\widetilde{V} \subset \widetilde{Y}$
as
in Lemma~\ref{lem:flatness1}.

For any $\tilde{y} \in \widetilde{V} (k)$,
let $\overline{\tilde{y}_{K}}$ be the geometric point 
of $Y$ corresponding to $\tilde{y}$
(cf. \S~\ref{subsec:relativeheight}).
Note that 
$\{ \tilde{y} \} \times \mathfrak{B}'
= \mathfrak{B}'$
naturally
and note that
$\overline{\tilde{y}_{K}}$
is also regarded as
the geometric generic point of $\{ \tilde{y} \} \times \mathfrak{B}'
= \mathfrak{B}'$
and is a $\overline{K}$-valued point.
Let
$f_{\tilde{y}} : \left\{ \tilde{y} \right\} \times \mathscr{A} 
\to 
\left\{ \tilde{y} \right\} \times \mathfrak{B}'$
be the restriction of 
$\pr_{\mathscr{B}}
: \mathscr{B} \times_{\mathfrak{B}'} \mathscr{A}
\to \mathscr{B}$.
We remark that, via natural identification 
$\left\{ \tilde{y} \right\} \times \mathscr{A} = \mathscr{A}$
and $\left\{ \tilde{y} \right\} \times \mathfrak{B}' = \mathfrak{B}'$,
$f_{\tilde{y}}$ coincides with $f$.
Further, we remark that the pair $\left( 
f_{\tilde{y}}
: \left\{ \tilde{y} \right\} \times \mathscr{A} 
\to 
\left\{ \tilde{y} \right\} \times
\mathfrak{B}' , 
\mathscr{L} \right)$
is a model of $\left( \left\{ \overline{\tilde{y}_{K}} \right\}
\times A
, L \right)$
as in Proposition~\ref{prop:model1},
where $\mathscr{L}$ and $L$ are 
regarded as line bundles on $\left\{ \tilde{y} \right\} \times \mathscr{A}$
and $\left\{ \overline{\tilde{y}_{K}} \right\}
\times A$ via the natural isomorphism
$\left\{ \tilde{y} \right\} \times \mathscr{A} = \mathscr{A}$
and
$\left\{ \overline{\tilde{y}_{K}} \right\}
\times A = A$ respectively.

Let $
q : \mathscr{B} = \widetilde{B} \times_{\Spec (k)} \mathfrak{B}'
\to \widetilde{B}$ be the natural projection.
For any $\tilde{y} \in \widetilde{V} (k)$,
set $\mathscr{X}_{\tilde{y}} := (
q \circ h)^{-1} (\tilde{y})$,
which
is a closed subscheme of $\{ \tilde{y} \} \times \mathscr{A}$.
Let $p : Y \times A \to Y$ be the restriction of $\pr_{B} : B \times A \to B$.
Then
we see that
the geometric fiber $X_{\overline{\tilde{y}_{K}}}$
of $p|_X : X \to Y$ over $\overline{\tilde{y}_{K}}$
equals
the geometric generic fiber of
$h |_{\mathscr{X}_{\tilde{y}}}
: \mathscr{X}_{\tilde{y}} \to
\{ \tilde{y} \} \times \mathfrak{B}'$.
By the condition (b) of Lemma~\ref{lem:flatness1},
$h |_{\mathscr{X}_{\tilde{y}}}
: \mathscr{X}_{\tilde{y}} \to
\{ \tilde{y} \} \times \mathfrak{B}' (= \mathfrak{B}')$
is flat over any point 
of $\mathfrak{B}'$
of codimension $1$.
Further, 
since $\mathscr{X}$ is an integral scheme
 and
since $\pr_{\widetilde{Y}} \circ h$ is flat over $\widetilde{V}$
by the condition (a) of Lemma~\ref{lem:flatness1},
$\mathscr{X}_{\tilde{y}}$ has pure dimension.
Set $d := \dim (X)$ and $e := \dim (Y)$.
Then, it follows from
Lemma~\ref{lem:height-intersection}
and the definition (\ref{eq:defofmathbfh})
of $\mathbf{h}^{L}_{X/Y}$
that,
for any $\tilde{y} \in \widetilde{V} (k)$,
we have
\addtocounter{Claim}{1}
\begin{align} \label{eq:height-inetersection-inproof}
[K':K]
\mathbf{h}^{L}_{X/Y} ( \tilde{y} )
=
\deg_{\mathcal{H}'}
\left(
f_{\tilde{y}}
\right)_{\ast}
\left(
\cherncl_{1}
(\mathscr{L})^{\cdot (d+1 - e)}
\cdot
\left[ 
\mathscr{X}_{\tilde{y}}
\right]
\right)
,
\end{align}
where $\mathcal{H}'$ is naturally regarded as a line bundle on
$\{ \tilde{y} \} \times \mathfrak{B}'$ via
$\{ \tilde{y} \} \times \mathfrak{B}' = \mathfrak{B}'$.
Equality (\ref{eq:height-inetersection-inproof})
will be used later.

Set
$M \boxtimes L := \pr_{B}^{\ast} (M) \otimes \pr_{A}^{\ast} (L)$,
where $\pr_{A} : B \times A \to A$ is the canonical projection.
The multilinearity of the canonical height with respect to
line bundles (cf. \cite[Theorem~11.18~(a)]{gubler0})
gives us
\begin{align*} 
\widehat{h}_{M \boxtimes L} (X)
=
\sum_{j = 0}^{d+1}
\left(
\begin{matrix}
d + 1 \\
j
\end{matrix}
\right)
\widehat{h}_{\tiny{\underbrace{\pr_{B}^{\ast} (M)
,\ldots , \pr_{B}^{\ast} (M)}_{j},
\underbrace{
\pr_{A}^{\ast} (L) , \ldots , \pr_{A}^{\ast} (L)}_{d+1-j}}} (X)
.
\end{align*}
Since $X$ has dense small points,
we have
$\widehat{h}_{M \boxtimes L} (X) = 0$ by Proposition~\ref{prop:dense-height0}.
On the other hand,
by \cite[Theorem~11.18~(e)]{gubler0},
we have
\[
\widehat{h}_{\tiny{\underbrace{\pr_{B}^{\ast} (M)
,\ldots , \pr_{B}^{\ast} (M)}_{j},
\underbrace{
\pr_{A}^{\ast} (L) , \ldots , \pr_{A}^{\ast} (L)}_{d+1-j}}} (X)
\geq 0
\]
for each $j$.
It follows that
\addtocounter{Claim}{1}
\begin{align} \label{eq:relativeheight1}
\widehat{h}_{\tiny{\underbrace{\pr_{B}^{\ast} (M)
,\ldots , \pr_{B}^{\ast} (M)}_{e},
\underbrace{
\pr_{A}^{\ast} (L) , \ldots , \pr_{A}^{\ast} (L)}_{d+1-e}}} (X)
= 0
.
\end{align}

Recall that
$\varphi
:
\mathscr{B} \times_{\mathfrak{B}'}
\mathscr{A} \to \mathfrak{B}'
$
is 
a model of
$B \times A$
and that
$q : \mathscr{B} = \widetilde{B} \times_{\Spec (k)} 
\mathfrak{B}' \to \widetilde{B}$
is the canonical projection.
Set 
$\mathscr{M} := q^{\ast} \left( \widetilde{M} \right)$.
Then      
$\left( 
\varphi
, \pr_{\mathscr{B}}^{\ast} (\mathscr{M})
\right)$
and 
$\left( 
\varphi
, \pr_{\mathscr{A}}^{\ast} (\mathscr{L})
\right)$
are models 
of
$( B \times A , \pr_{B}^{\ast} (M))$
and
$( B \times A , \pr_A^{\ast} (L))$
respectively,
satisfying the 
conditions in Proposition~\ref{prop:model1}
(cf. Remark~\ref{rem:modelofconstant}).
Via 
the description (\ref{eq:canonicalheight})
or Lemma~\ref{lem:height-intersection},
equality
(\ref{eq:relativeheight1})
 gives us
\addtocounter{Claim}{1}
\begin{align} \label{eq:relativeheight2}
\deg_{\mathcal{H}'}
\varphi_{\ast}
\left(
\cherncl_{1}
\left(\pr_{\mathscr{A}}^{\ast}(\mathscr{L})\right)^{\cdot (d+1 - e)}
\cdot
\cherncl_{1}
\left( \pr_{\mathscr{B} }^{\ast} (\mathscr{M}) \right)^{\cdot e}
\cdot
[\mathscr{X}]
\right)
=
0
.
\end{align}

With the preparation so far, let us prove 
that $\mathbf{h}^{L}_{X/Y}$ vanishes
generically on $\widetilde{Y} (k)$.
Let $\left| \widetilde{M} \right|$ denote the
complete linear system over $\widetilde{B}$
associated to $\widetilde{M}$.
Since $\widetilde{M}$ is very ample,
the linear system $(q \circ \pr_{\mathscr{B}})^{\ast} \left| \widetilde{M} \right| $
over $\mathscr{B} \times_{\mathfrak{B}'} \mathscr{A}$
is base-point free.
By Lemma~\ref{lem:representingcycle},
it follows that,
for a general $( \widetilde{D}_1 , \ldots , \widetilde{D}_e )
\in \left| \widetilde{M} \right|^{e} $,
the cycle
\[
\left[
(q \circ \pr_{\mathscr{B}})^{-1} 
\left( 
\widetilde{D}_1
\right)
\cap 
\cdots
\cap
(q \circ \pr_{\mathscr{B}})^{-1} 
\left( 
\widetilde{D}_e
\right)
\cap
\mathscr{X}
\right]
=
\left[
(q \circ \pr_{\mathscr{B}})^{-1} 
\left(
\widetilde{D}_1
\cap 
\cdots
\cap
\widetilde{D}_e
\right)
\cap
\mathscr{X}
\right]
\]
represents the cycle class
$
\cherncl_1 \left( 
\pr_{\mathscr{B}}^{\ast} (\mathscr{M})
\right)^{\cdot e}
\cdot
[\mathscr{X}]
$.
Furthermore,
for a general $( \widetilde{D}_1 , \ldots , \widetilde{D}_e )
\in \left| \widetilde{M} \right|^{e} $,
the intersection
$\widetilde{D}_1
\cap 
\cdots
\cap
\widetilde{D}_e$
is a finite 
closed subscheme
of the regular locus
$\widetilde{V}_{\mathrm{reg}}$ of $\widetilde{V} \subset \widetilde{Y}$,
so that
we write
\[
\left[
\widetilde{D}_1
\cap 
\cdots
\cap
\widetilde{D}_e
\right]
=
\left[ \tilde{y}_1 \right]
+ \cdots + 
\left[ \tilde{y}_m
\right]
\]
with $\tilde{y}_1 , \ldots , \tilde{y}_m
\in \widetilde{V}_{\mathrm{reg}} (k)$.
Then we have
\[
\left[
(q \circ \pr_{\mathscr{B}})^{-1} 
\left(
\widetilde{D}_1
\cap 
\cdots
\cap
\widetilde{D}_e
\right)
\cap
\mathscr{X}
\right]
=
\sum_{i = 1}^{m}
\left[
\left( \{ \tilde{y}_i \} \times \mathscr{A} \right)
\cap \mathscr{X}
\right]
,
\]
and hence
\addtocounter{Claim}{1}
\begin{multline} \label{eq:heightsum}
\deg_{\mathcal{H}'}
\varphi_{\ast}
\left(
\cherncl_{1}
\left(
\pr_{\mathscr{A}}^{\ast} (\mathscr{L})
\right)^{\cdot (d+1 - e)}
\cdot
\cherncl_{1}
\left(
\pr_{\mathscr{B}}^{\ast}
(\mathscr{M})
\right)^{\cdot e}
\cdot
[\mathscr{X}]
\right)
\\=
\sum_{i = 1}^{m}
\deg_{\mathcal{H}'}
\left(
\varphi|_{\{ \tilde{y}_i \} \times \mathscr{A}} 
\right)_{\ast}
\left(
\cherncl_{1}
(\mathscr{L})^{\cdot (d+1 - e)}
\cdot
\left[
\left( \{ \tilde{y}_i \} \times \mathscr{A} \right)
\cap \mathscr{X}
\right]
\right)
,
\end{multline}
where $\mathscr{L}$ is naturally regarded as a line bundle on
$\{ \tilde{y}_i \} \times \mathscr{A}$.
%

Since
$( q \circ \pr_{\mathscr{B}}) |_{\mathscr{X}} = 
q \circ h$
by the definition of $h$,
we have
$\left( 
\{ \tilde{y}_i \} \times \mathscr{A} \right) \cap \mathscr{X}
= \mathscr{X}_{\tilde{y}_i}$.
We remark that 
$f_{\tilde{y}_i} = \varphi|_{\{ \tilde{y}_i \} \times \mathscr{A}} $
via the natural isomorphism 
$\{ \tilde{y}_i \} \times \mathfrak{B}' = \mathfrak{B}'$.
By (\ref{eq:height-inetersection-inproof}),
it follows that 
the
right-hand side of (\ref{eq:heightsum})
 equals $
[K' :K]\sum_{i = 1}^{m}
\mathbf{h}^{L}_{X/Y} (\tilde{y}_{i})$.
On the other hand,
the left-hand side
of (\ref{eq:heightsum}) equals $0$
by 
(\ref{eq:relativeheight2}).
By
Remark~~\ref{rem:mathbfh-nonnegative},
we thus have
$
\mathbf{h}^{L}_{X/Y} (\tilde{y}_{i}) = 0
$.

In summary, we have seen that
there exists a dense open subset $\mathcal{U} 
\subset
\left| \widetilde{M} \right|^{e}$ such that,
for any
$
(\widetilde{D}_{1} , \ldots ,\widetilde{D}_e)
\in \mathcal{U}$,
and
for any
$\tilde{y} \in \widetilde{D}_1
\cap 
\cdots
\cap
\widetilde{D}_e$
we have
\addtocounter{Claim}{1}
\begin{align} \label{eq:heightofyi}
\mathbf{h}^{L}_{X/Y} (\tilde{y}) = 0.
\end{align}
Since 
$\left| \widetilde{M} \right|$ is very ample,
a general point $\tilde{y} \in \widetilde{Y} (k)$
is a point of the intersection
$\widetilde{D}_1
\cap 
\cdots
\cap
\widetilde{D}_e
$
for some
$(\widetilde{D}_{1} , \ldots ,\widetilde{D}_e)
\in \mathcal{U}$.
It follows from
equality
(\ref{eq:heightofyi})
that
$\mathbf{h}^{L}_{X/Y} (\tilde{y}) = 0$ for a general $\tilde{y}
\in \widetilde{Y} (k)$.
This is the second statement of the proposition,
and thus
we complete the proof.
\QED

\begin{Proposition} \label{prop:final2}
Let $A$ be an abelian variety over $\overline{K}$
with trivial $\overline{K}/k$-trace
and let $L$ be an even ample line bundle on $A$.
Let $\widetilde{Y}$ be an irreducible variety over $k$
and set
$Y := \widetilde{Y} \otimes_{k} \overline{K}$.
Let
$p : Y \times A \to Y$ denote the canonical projection
and
let 
$X$ be an irreducible closed subvariety of $Y \times A$
such that
the restriction 
$p|_{X} : X \to Y$ 
is surjective.
Suppose that there exists a dense open subset $\widetilde{U}$
of $\tilde{Y}$ with
$\widetilde{U} \subset \widetilde{Y}_{\mathrm{pd}}$
such that $\mathbf{h}^{L}_{X/Y} = 0$ over $\widetilde{U}$.
Then, if the geometric Bogomolov conjecture holds for $A$,
then
there exists a torsion subvariety $T$ of $A$
such that $X = Y \times T$.
\end{Proposition}

\Proof
To use Proposition~\ref{prop:key1},
we 
check that the assumptions of it
are satisfied.
First, since the generic point $\tilde{\eta}$ of $\widetilde{Y}$
sits in $\widetilde{U}$,
we have $\mathbf{h}^{L}_{X/Y} (\tilde{\eta}) = 0$.
Thus the condition (a) of Proposition~\ref{prop:key1} is verified.
To see the condition (b),
let $S$ be the preimage of $\widetilde{U} (k)$ by the natural morphism
$Y \to \widetilde{Y}$.
Note that $S$ is a dense subset of $Y$.
For any $\tilde{y} \in \widetilde{U}(k)$,
let $\overline{\tilde{y}_K}$ be the corresponding geometric 
point of $Y$.
By
definition,
we have
$\overline{\tilde{y}_K} \in S \cap Y \left( \overline{K} \right)$.
Since $\mathbf{h}^{L}_{X/Y} (\tilde{y}) = 0$,
it follows from Remark~\ref{rem:mathbfh-nonnegative}
that
any irreducible component of $X_{\overline{\tilde{y}_K}}$
is a closed subvariety of $A$ of
canonical height zero.
Suppose that the geometric Bogomolov conjecture holds
for $\{ \tilde{y} \}
\times A = A$.
Then, it follows from Proposition~\ref{prop:dense-height0}
that
$X_{\overline{\tilde{y}_K}}$ is a torsion subvariety.
Thus, the condition (b) of Proposition~\ref{prop:key1}
is verified,
and hence
Proposition~\ref{prop:final2} 
follows from Proposition~\ref{prop:key1}.
\QED

As a consequence, we obtain the following theorem.

\begin{Theorem} \label{thm:height-torsion}
Let $A$ be a nowhere degenerate abelian variety over $\overline{K}$
with
trivial $\overline{K}/k$-trace.
Let
$\widetilde{B}$ be an abelian variety over $k$
and set
$B = \widetilde{B} \otimes_k \overline{K}$.
Let $\pr_A : B \times A \to A$ and $\pr_B : B \times A \to B$ 
denote the canonical projections.
Let $X$ be an irreducible closed subvariety of $B \times A$
and set $Y := \pr_B (X)$ and $T:= \pr_A (X)$.
Suppose that $X$ has dense small points and assume that
the geometric Bogomolov conjecture holds for $A$.
Then,
$Y$ is a constant subvariety of $B$,
$T$ is a torsion subvariety of $A$,
and $X = Y \times T$ holds.
\end{Theorem}

\Proof
Let $L$ be an even ample line bundle on $A$.
It follows from Proposition~\ref{prop:final}
that $Y$ is a constant subvariety
and 
$\mathbf{h}^{L}_{X/Y} (y) = 0$ holds for
general $y \in \widetilde{Y} (k)$,
where $\widetilde{Y}$ is a closed subvariety of $\widetilde{B}$ with
$Y = \widetilde{Y} \otimes_k \overline{K}$.
By Proposition~\ref{prop:hgenericallyconstant},
therefore,
$\mathbf{h}^{L}_{X/Y} = 0$ on some dense open subset
of $\widetilde{Y}$ contained in $\widetilde{Y}_{\mathrm{pd}}$.
Then
Proposition~\ref{prop:final2}
concludes
that $T$ is a torsion subvariety,
and $X = Y \times T$ holds.
\QED

Let $\mathfrak{m}$ be the maximal nowhere degenerate abelian 
subvariety of $A$.
Let $\left( \widetilde{A}^{\overline{K}/k} , \Tr_{A} \right)$
be the $\overline{K}/k$-trace of $A$
and let $\mathfrak{t}$ be the image of $\Tr_{A}$.
Then,
by \cite[Proposition~7.11]{yamaki6},
we have $\mathfrak{t} \subset \mathfrak{m}$.

\begin{Remark} \label{rem:quotient-trivial-trace}
The quotient $\mathfrak{m} / \mathfrak{t}$
is nowhere degenerate and has trivial $\overline{K}/k$-trace.
Indeed, the nondegeneracy follows from
\cite[Lemma~7.8~(2)]{yamaki6},
and the triviality of the trace
follows from the well-known fact
that a surjective homomorphism between
abelian varieties over $\overline{K}$
induces a surjective homomorphism between their $\overline{K}/k$-traces
(cf. \cite[Lemma~1.5]{yamaki5}).
\end{Remark}

Now, we establish the second main theorem of this paper,
which 
includes
Theorem~\ref{thm:reduction-GBCforMNDAS}.

\begin{Theorem} [Theorem~\ref{thm:main2intro}] \label{thm:main2}
Let $A$ be an abelian variety over $\overline{K}$.
Let $\mathfrak{m}$ be the maximal nowhere degenerate abelian subvariety
of $A$ and let $\mathfrak{t}$ be the image of the 
$\overline{K}/k$-trace homomorphism
of $A$.
Then the following statements are equivalent to each other.
\begin{enumerate}
\renewcommand{\labelenumi}{(\alph{enumi})}
\item
The geometric Bogomolov conjecture holds for $A$.
\item
The geometric Bogomolov conjecture holds for $\mathfrak{m}$.
\item
The geometric Bogomolov conjecture holds for
$\mathfrak{m} / \mathfrak{t}$.
\end{enumerate}
\end{Theorem}

\Proof
The equivalence between (a) and (b) is nothing but 
Theorem~\ref{thm:reduction-GBCforMNDAS}.
It follows from \cite[Lemma~7.7]{yamaki6}
that (b) implies (c).

To prove that (c) implies (b),
suppose that
the
geometric Bogomolov conjecture holds for
$\mathfrak{m} / \mathfrak{t}$.
We put $B := \widetilde{A}^{\overline{K}/k} \otimes_{k} \overline{K}$.
Since $\mathfrak{t} \subset \mathfrak{m}$ 
and since $\Tr_A :
B \to \mathfrak{t}$
is an isogeny,
$\mathfrak{m}$ is isogenous to 
$B \times 
\mathfrak{m} / \mathfrak{t}$.
By \cite[Corollary~7.6]{yamaki6},
it suffices to
show that the conjecture holds for
this abelian variety.
Let $X \subset B \times 
\mathfrak{m} / \mathfrak{t}$
be an irreducible closed subvariety
having dense small points.
Let $\pr_B : B \times \mathfrak{m} / \mathfrak{t}
\to B$ and $\pr_{\mathfrak{m} / \mathfrak{t}} : B \times \mathfrak{m} / \mathfrak{t}
\to \mathfrak{m} / \mathfrak{t}$ be the natural projections.
Since $\mathfrak{m} / \mathfrak{t}$ 
is nowhere degenerate and
has trivial $\overline{K}/k$-trace 
(cf. Remark~\ref{rem:quotient-trivial-trace}),
Theorem~\ref{thm:height-torsion} then
tells us that $\pr_{B} (X)$ is a constant subvariety,
$\pr_{\mathfrak{m} / \mathfrak{t}} (X)$ is a torsion subvariety,
and $X = \pr_{B} (X) \times \pr_{\mathfrak{m} / \mathfrak{t}} (X)$.
This shows that the geometric Bogomolov conjecture holds for
$B \times 
\mathfrak{m} / \mathfrak{t}$.
Thus we 
conclude that (c) implies (b), which completes the proof
the theorem.
\QED

\begin{Remark} \label{rem:1from2}
\begin{enumerate}
\item
Theorem~\ref{thm:main1intro} follows from
Theorem~\ref{thm:main2}
(Theorem~\ref{thm:main2intro}).
Indeed, since the trace homomorphism is an isogeny,
$\dim \left( \widetilde{A}^{\overline{K}/k} \right) = \ndr (A)$
implies $\mathfrak{m} / \mathfrak{t} = 0$.
Then the geometric Bogomolov conjecture for $A$ holds by 
Theorem~\ref{thm:main2}
(Theorem~\ref{thm:main2intro}).
\item
Since the geometric Bogomolov conjecture holds for abelian varieties
of dimension not greater than $1$,
a direct application of Theorem~\ref{thm:main2}
tells  us 
that the conjecture holds for $A$ with
$\dim ( \mathfrak{m} / \mathfrak{t}) \leq 1$, or equivalently,
with
$\dim \left( \widetilde{A}^{\overline{K}/k} \right) \geq \ndr (A) -1$.
This
seems to give us a result stronger than 
Theorem~\ref{thm:main1intro}.
However,
it is not stronger in fact
because
$\dim \left( \widetilde{A}^{\overline{K}/k} \right) \geq \ndr (A) -1$
leads us to
$\dim \left( \widetilde{A}^{\overline{K}/k} \right) = \ndr (A)$.
Indeed, suppose
$\dim \left( \widetilde{A}^{\overline{K}/k} \right) \geq \ndr (A) -1$,
i.e.,
$\dim (\mathfrak{m} / \mathfrak{t}) \leq 1$.
If we had $\dim (\mathfrak{m} / \mathfrak{t}) = 1$,
then $\mathfrak{m} / \mathfrak{t}$ would be a constant variety
(see the argument in Remark~\ref{rem:genthmD}),
but this is a contradiction by Remark~\ref{rem:quotient-trivial-trace}.
Thus 
$\dim \left( \widetilde{A}^{\overline{K}/k} \right) \geq \ndr (A) -1$
implies $\dim (\mathfrak{m} / \mathfrak{t}) = 0$,
i.e.,
$\dim \left( \widetilde{A}^{\overline{K}/k} \right) = \ndr (A)$.
\end{enumerate}
\end{Remark}

\begin{Remark} \label{rem:reductiontonew1}
Since
the abelian variety $\mathfrak{m} / \mathfrak{t}$
is a nowhere degenerate abelian variety over $\overline{K}$
with trivial $\overline{K}/k$-trace
(cf. Remark~\ref{rem:quotient-trivial-trace}),
it follows from Theorem~\ref{thm:main2}
that
Conjecture~\ref{GBCforAV} is reduced to 
the geometric Bogomolov conjecture for nowhere degenerate
abelian varieties with trivial $\overline{K}/k$-trace.
Further,
since
any special subvariety of an abelian variety
with trivial $\overline{K}/k$-trace is a torsion subvariety,
Conjecture~\ref{GBCforAV}
is in fact reduced to the following conjecture\footnote{This should be compared
with  \cite[Conjecture~7.22]{yamaki6}.}.
\end{Remark}

\begin{Conjecture} [Geometric Bogomolov conjecture
for nowhere abelian varieties with trivial trace] \label{conj:new1}
Let $A$ be a nowhere degenerate abelian variety over $\overline{K}$
with trivial $\overline{K}/k$-trace.
Then
any irreducible closed subvariety of $A$
with
dense small points
is a torsion subvariety.
\end{Conjecture}

\renewcommand{\thesection}{Appendix} 
\renewcommand{\theTheorem}{A.\arabic{Theorem}}
\renewcommand{\theClaim}{A.\arabic{Theorem}.\arabic{Claim}}
\renewcommand{\theequation}{A.\arabic{Theorem}.\arabic{Claim}}
\renewcommand{\theProposition}
{A.\arabic{Theorem}.\arabic{Proposition}}
\renewcommand{\theLemma}{A.\arabic{Theorem}.\arabic{Lemma}}
\setcounter{section}{0}
\renewcommand{\thesubsection}{A.\arabic{subsection}}


\section{Field extension and the trace}

In this appendix,
we show the following lemma,
which is used in the proof of Proposition~\ref{prop:key1}.

\begin{Lemma} \label{lem:trivial-trivial}
Let $F/k$ 
be a field extension
%
with $F$ 
algebraically closed
and
let 
$\mathfrak{k}/k$ be a field extension 
with $\mathfrak{k}$ algebraically closed
such that $\mathfrak{k} \otimes_{k} F$ is 
an integral domain.
Let $\mathfrak{F}$ be a field containing 
$\mathfrak{k} \otimes_{k} F$ as a subring.
Let $A$ be an abelian variety over $F$
and suppose that $A$ has trivial $F/k$-trace.
Then
$A \times_{\Spec (F)} \Spec (\mathfrak{F})$ has trivial 
$\mathfrak{F}/\mathfrak{k}$-trace.
\end{Lemma}

\Proof
Let $B$ be an abelian variety over $\mathfrak{k}$.
Then there
exist a finitely generated $k$-subalgebra $R \subset \mathfrak{k}$ 
and an abelian scheme $\mathcal{B} \to \Spec (R)$
such that
$B = \mathcal{B} \times_{\Spec (R)} \Spec (\mathfrak{k})$.
Let 
\[
\phi : B \times_{\Spec (\mathfrak{k})} \Spec (\mathfrak{F}) \to A
\times_{\Spec (F)} \Spec (\mathfrak{F})
\] 
be a homomorphism.
Then
there exist a finitely generated $F$-algebra $S$
with
$R  \otimes_k F \subset S \subset
\mathfrak{F}$ 
and
a homomorphism 
\[
\Phi : \mathcal{B} \times_{\Spec (R)} \Spec (S)
\to A \times_{\Spec \left( F \right)} \Spec (S)
\]
such that its base-change 
$\Phi_{\mathfrak{F}}$ to $\Spec ( \mathfrak{F})$
coincides with $\phi$.

We set $X' := \Spec (R)$, $X := \Spec (R \otimes_k F)
= X' \times_{\Spec (k)} \Spec \left( F \right)$,
and $Y := \Spec (S)$.
The morphism $f : Y \to X$ induced from the inclusion
$R \otimes_k F \subset S$ is a morphism of varieties over 
$F$
and is dominant,
so that there exists
a dense open subset $U \subset X$ which is contained
in the image of $f$.
Since $X
= X' \times_{\Spec (k)} \Spec \left( F \right)$,
we have a natural injection $X' (k) \to X \left( 
F \right)$,
and
let $X(k)$ denote its image.
Then $f^{-1} \left( X(k) \cap U \right)$ is dense in $Y$.

Let $g : X \to X'$ be the natural projection.
For any $x \in X(k) \cap U$,
the fiber $f^{-1}(x)$ is a closed subscheme of $Y$
and is a scheme over $X'$ via $g \circ f$.
We consider the restriction
\[
\Phi_{f^{-1} (x)} : \mathcal{B} \times_{X'} f^{-1} (x) \to 
A \times_{\Spec \left( F \right)} f^{-1} (x)
\] 
of $\Phi$,
which
is a homomorphism of abelian schemes over $f^{-1} (x)$.
Let
$y \in f^{-1} (x)$ be any closed point,
which is a $F$-valued point of $Y$.
Then
the fiber $\mathcal{B} \times_{X'} \{ y \}$ of the abelian scheme
$\mathcal{B} \times_{X'} f^{-1} (x) \to f^{-1} (x)$ coincides with 
\[
\mathcal{B} \times_{X'} \{ g (x)\} \times_{\{ g (x) \}} \{ y \}
= \left( \mathcal{B} \times_{X'} \{ g (x)\} \right) \otimes_{k} 
F.
\]
Further,
the fiber of $A \times_{\Spec \left( F \right)} f^{-1} (x) 
\to f^{-1} (x) $ over $y$
equals $A$.
Since $\mathcal{B} \times_{X'} \{ g (x)\}$ is an abelian variety over $k$
and since $A$ has trivial $F / k$-trace,
it follows that
the homomorphism
$\Phi_{f^{-1} (x)}$
is trivial over $y$.
This means that
$\Phi_{f^{-1} (x)}$ is trivial over any closed point of $f^{-1} (x)$,
and hence
$\Phi_{f^{-1} (x)}$ itself is the trivial homomorphism.
Since $x$ is any point of $ X(k) \cap U$
and
$f^{-1} \left( X(k) \cap U \right)$ is dense in $Y$,
it follows further that $\Phi$ is the trivial homomorphism.
Thus $\phi$, which is the base-change of $\Phi$ to $\mathfrak{F}$,
is also trivial.
\QED

Let $F$, $\mathfrak{k}$
and
$\mathfrak{F}$
be as in Lemma~\ref{lem:trivial-trivial}.
Let $\left( A^{F/k} , \Tr_A \right)$ be the $F/k$-trace of $A$
and let 
$\left( A_{\mathfrak{F}}^{\mathfrak{F}/\mathfrak{k}} , \Tr_{A_{\mathfrak{F}}} \right)$ be the 
$\mathfrak{F}/\mathfrak{k}$-trace 
of $A_{\mathfrak{F}} := A \otimes_F \mathfrak{F}$.
Let $\phi : A^{F/k} \otimes_{k} \mathfrak{k}
\to A_{\mathfrak{F}}^{\mathfrak{F}/\mathfrak{k}}$ 
be the homomorphism induced by the universality from the base-change
$\Tr_A \otimes_{F} \mathfrak{F} : 
A^{F/k} \otimes_{k} \mathfrak{F} \to A_{\mathfrak{F}}$ of $\Tr_A$.
We end with a remark that
$\phi$ is a purely inseparable isogeny.
Since the trace homomorphism is purely inseparable to its image
(cf. \cite[VIII \S~3 Corollary~2]{lang1}
or
\cite[Lemma~1.4]{yamaki5}),
the same holds for $\phi$. Therefore,
we only have to show that 
$\dim \left( A^{F/k} \otimes_{k} \mathfrak{k} \right)
\geq
\dim \left( A_{\mathfrak{F}}^{\mathfrak{F}/\mathfrak{k}} \right)$.
Let 
$\mathfrak{t} \subset A$ be the image of $\Tr_A$
and let $q : A \to  A / \mathfrak{t}$ be the quotient.
By the same argument as Remark~\ref{rem:quotient-trivial-trace},
we note that
$A / \mathfrak{t}$ has trivial $F / k$-trace.
Thus
Lemma~\ref{lem:trivial-trivial}
tells us that
$(A / \mathfrak{t}) \otimes_{F} \mathfrak{F}$ has
trivial $\mathfrak{F} / \mathfrak{k}$-trace.
It follows that the composite
\[
\begin{CD}
A_{\mathfrak{F}}^{\mathfrak{F}/\mathfrak{k}} \otimes_{\mathfrak{k}} \mathfrak{F}
@>{\Tr_{A_{\mathfrak{F}}}}>> A_{\mathfrak{F}} 
@>{q \otimes_{\mathfrak{k}} {\mathfrak{F}}}>> 
(A / \mathfrak{t}) \otimes_{F} \mathfrak{F}
\end{CD}
\]
is trivial,
which means that 
$\Tr_{A_{\mathfrak{F}}} \left( 
A_{\mathfrak{F}}^{\mathfrak{F}/\mathfrak{k}} \otimes_{\mathfrak{k}} \mathfrak{F}
\right) \subset \mathfrak{t} \otimes_{k} \mathfrak{k}$.
Since
the trace homomorphisms
are finite
(cf. \cite[Lemma~1.4]{yamaki5}), we obtain
\[
\dim \left( A^{F/k} \otimes_{k} \mathfrak{k} \right)
=
\dim \left( \mathfrak{t} \otimes_{k} \mathfrak{k} \right)
\geq
\dim \left(
\Tr_{A_{\mathfrak{F}}} \left( 
A_{\mathfrak{F}}^{\mathfrak{F}/\mathfrak{k}} \otimes_{\mathfrak{k}} \mathfrak{F}
\right)
\right)
=
\dim
\left( 
A_{\mathfrak{F}}^{\mathfrak{F}/\mathfrak{k}} \otimes_{\mathfrak{k}} \mathfrak{F}
\right)
=
\dim
\left( A_{\mathfrak{F}}^{\mathfrak{F}/\mathfrak{k}} \right)
,
\]
as required.


\small{

}

\end{document}